\input amstex
\documentstyle{amsppt} 
\magnification=\magstep1
\NoRunningHeads
\NoBlackBoxes
\topmatter
\title Orbit equivalence of one-sided subshifts 
and the associated $C^*$-algebras
\endtitle
\author Kengo Matsumoto
\endauthor
\affil Department of Mathematical Sciences \\
Yokohama City University \\
Seto 22-2, Kanazawa-ku, Yokohama 236-0027, JAPAN
\endaffil
\abstract{A $\lambda$-graph system ${\frak L}$
is a generalization of a finite 
labeled graph and presents a subshift. 
We will prove that the topological dynamical systems 
$(X_{{\frak L}_1},\sigma_{{\frak L}_1})$ and 
$(X_{{\frak L}_2},\sigma_{{\frak L}_2})$ for $\lambda$-graph systems 
${\frak L}_1$ and ${\frak L}_2$
are continuously orbit equivalent 
if and only 
if there exists an isomorphism between 
the associated $C^*$-algebras 
${\Cal O}_{{\frak L}_1}$ and 
${\Cal O}_{{\frak L}_2}$ keeping their commutative 
$C^*$-subalgebras 
$C(X_{{\frak L}_1})$ and $C(X_{{\frak L}_2})$.
It is also equivalent to the condition that there exists a homeomorphism
from
$X_{{\frak L}_1}$ to $X_{{\frak L}_2}$ 
intertwining their topological full inverse semigroups.  
In particular,
one-sided subshifts 
$X_{\Lambda_1}$ and $X_{\Lambda_2}$
are $\lambda$-continuously orbit equivalent 
if and only 
if there exists an isomorphism between 
the associated $C^*$-algebras ${\Cal O}_{\Lambda_1}$ 
and 
${\Cal O}_{\Lambda_2}$ keeping their commutative 
$C^*$-subalgebras $C(X_{\Lambda_1})$ and $C(X_{\Lambda_2})$.
}
\endabstract
\endtopmatter


\def\Zp{{ {\Bbb Z}_+ }}

\def\OL{{{\Cal O}_{\frak L}}}
\def\OLA{{{\Cal O}_{{\frak L}_1}}}
\def\OLB{{{\Cal O}_{{\frak L}_2}}}

\def\DL{{{\Cal D}_{\frak L}}}
\def\DLA{{{\Cal D}_{{\frak L}_1}}}
\def\DLB{{{\Cal D}_{{\frak L}_2}}}

\def\FL{{{\Cal F}_{\frak L}}}

\def\OLambda{{{\Cal O}_{\Lambda}}}
\def\OLambdaA{{{\Cal O}_{\Lambda_1}}}
\def\OLambdaB{{{\Cal O}_{\Lambda_2}}}

\def\DLambda{{{\frak D}_\Lambda}}
\def\DLambdaA{{{\frak D}_{\Lambda_1}}}
\def\DLambdaB{{{\frak D}_{\Lambda_2}}}

\def\FKL{{ {\Cal F}_k^{l} }}

\def\OA{{{\Cal O}_A}}
\def\OB{{{\Cal O}_B}}

\def\DA{{{\frak D}_A}}
\def\DB{{{\frak D}_B}}

\def\Homeo{{{\operatorname{Homeo}}}}

\def\Aut{{{\operatorname{Aut}}}}

\def\id{{{\operatorname{id}}}}
\def\supp{{{\operatorname{supp}}}}


\heading 1. Introduction
\endheading
H. Dye has initiated to study of orbit equivalence of ergodic finite measure  preserving transformations, who proved that any two such transformations are orbit equivalent ([D], [D2]).
W. Krieger [Kr] has proved that two ergodic non-singular transformations are orbit equivalent if and only if the associated von Neumann crossed produtcs are isomorphic.
 In topological setting, Giordano-Putnam-Skau [GPS],[GPS2] (cf.[HPS]) have proved that two Cantor minimal systems are strong orbit equivalent 
 if and only if the associated $C^*$-crossed products are isomorphic.
In more general setting, 
J. Tomiyama [To]  (cf. [BT], [To2] ) has proved that    
two topological free homeomorphisms
$(X,\phi)$ and $(Y,\psi)$ on compact Hausdorff spaces
are continuously orbit equivalent if and only if 
there exists an isomorphism between the associated  $C^*$-crossed products 
keeping their commutative $C^*$-subalgebras $C(X)$ and $C(Y)$.
He also proved that it is equivalent to the condition that 
there exists a homeomorphism $h : X \rightarrow Y$ such that
$ h $ preserves their topological full groups.
Orbit equivalence of continuous maps on compact Hausdorff spaces 
that are not homeomorphisms are not covered by the above Tomiyama's setting.
The class of one-sided subshifts
 is an important class of topological dynamical systems on Cantor sets
with continuous surjections
that are not homeomorphisms.
 The one-sided topological Markov shifts is a subclass of the class.
 The associated $C^*$-algebras to the topological Markov shifts
are known to be  the Cuntz-Krieger algebras.
In the recent paper [Ma5], 
the author has shown that  
similar results to the Tomiyama's results hold
for one-sided topological Markov shifts.
He has proved that 
 one-sided topological Markov shifts
$(X_A,\sigma_A)$ and $(X_B,\sigma_B)$ for matrices $A$ and $B$ with 
entries in $\{0,1\}$ 
are continuously orbit equivalent 
if and only 
if there exists an isomorphism between 
the Cuntz-Krieger algebras ${\Cal O}_A$ and ${\Cal O}_B$ 
keeping their commutative 
$C^*$-subalgebras $C(X_A)$ and $C(X_B)$
( Note that the term \lq\lq topological \rq\rq orbit equivalence has been used in [Ma5] instead of \lq\lq continuous \rq\rq orbit equivalence).
It is also equivalent to the condition that there exists a homeomorphism
from
$X_A$ to $X_B$ intertwining their topological full groups
$[\sigma_A]_c$ and $[\sigma_B]_c$.

In this paper 
we will extend the above results for one-sided topological Markov shifts to the class of general one-sided subshifts.
A $\lambda$-graph system ${\frak L}$
is a generalization of a finite 
labeled graph and presents a subshift.
It yields a topological dynamical system 
$(X_{\frak L},\sigma_{\frak L})$
of a zero-dimensional compact Hausdorff space $X_{\frak L}$ with 
shift transformation $\sigma_{\frak L}$, 
that is a continuous surjection and not a homeomorphism. 
The $C^*$-algebra $\OL$ is associated with the dynamical system 
$(X_{\frak L},\sigma_{\frak L})$
such that $C(X_{\frak L})$ is naturally embedded into $\OL$ 
as a diagonal algebra of the canonical AF-algebra $\FL$ inside of 
$\OL$. 
We will prove that the topological dynamical systems 
$(X_{{\frak L}_1},\sigma_{{\frak L}_1})$ and 
$(X_{{\frak L}_2},\sigma_{{\frak L}_2})$ for $\lambda$-graph systems 
${\frak L}_1$ and ${\frak L}_2$
are continuously orbit equivalent 
if and only 
if there exists an isomorphism between 
the associated $C^*$-algebras 
${\Cal O}_{{\frak L}_1}$ and 
${\Cal O}_{{\frak L}_2}$ keeping their commutative 
$C^*$-subalgebras $C(X_{{\frak L}_1})$ and $C(X_{{\frak L}_2})$.
It is also equivalent to the condition that there exists a homeomorphism
from
$X_{{\frak L}_1}$ to $X_{{\frak L}_2}$ 
intertwining their topological full inverse semigroups
$[\sigma_{{\frak L}_1}]_{sc}$ 
and 
$[\sigma_{{\frak L}_1}]_{sc}$.  
Let 
$X_{\Lambda_1}$ and $ X_{\Lambda_2}$
be the right one-sided subshifts for two-sided subshifts
$\Lambda_1$ and $\Lambda_2$ respectively.
We in particular show that two one-sided subshifts 
$X_{\Lambda_1}$ and $X_{\Lambda_2}$
are $\lambda$-continuously orbit equivalent 
if and only 
if there exists an isomorphism between 
the associated $C^*$-algebras ${\Cal O}_{\Lambda_1}$ 
and 
${\Cal O}_{\Lambda_2}$ keeping their commutative 
$C^*$-subalgebras $C(X_{\Lambda_1})$ and $C(X_{\Lambda_2})$,
where
${\Cal O}_{\Lambda_1}$ 
and 
${\Cal O}_{\Lambda_2}$ are the $C^*$-algebras associated with subshifts ([Ma], cf. [CaM]).

Let $[\sigma_{\frak L}]_c$ be the topological full group of 
$(X_{\frak L},\sigma_{\frak L})$ 
whose elements consist of homeomorphisms $\tau$ on $X_{\frak L}$ 
such that $\tau (x) $ 
is contained in the orbit $orb_{\sigma_{\frak L}}(x)$ of $x$
under $\sigma_{\frak L}$ for $x \in X_A$, 
and its orbit cocycles are continuous. 
If ${\frak L}$ comes from a finite directed graph and hence 
$X_{\frak L}$ is a topological Markov shift,
then the topological full group 
is large enough to cover orbits of $x \in X_{\frak L}$.
However if $\frak L$ does not come from a finite graph,
the topological full group is not necessarily large enough
to cover orbits of $X_{\frak L}$.
To obtain enough informations of orbit structure 
of
$(X_{\frak L},\sigma_{\frak L})$,
we need to enlarge $[\sigma_{\frak L}]_c$ 
to
topological inverse semigroup $[\sigma_{\frak L}]_{sc}$
whose elements consist of partial homeomorphisms
$\tau$ on $X_{\frak L}$ 
such that 
$\tau(x)$ 
is contained in $orb_{\sigma_{\frak L}}(x)$
for each $x$ in the domain of $\tau$.
Let us denote by $\DL$ the commutative $C^*$-subalgebra
$C(X_{\frak L})$ of $\OL$. 
The corresponding object to the inverse semigroup 
$[\sigma_{\frak L}]_{sc}$ is the normalizer semigroup
$N_s(\OL,\DL)$ of $\DL$ in $\OL$ whose elements consist of
partial isometries $v $ of $\OL$ such that
$v \DL v^* \subset \DL$
and
$v^* \DL v \subset \DL$.
Then we will show that the exact sequence 
$$
1 
\longrightarrow
{\Cal U}(\DL)
\longrightarrow
N_s(\OL,\DL)
\longrightarrow
[\sigma_{\frak L}]_{sc}
\longrightarrow
1
$$
of semigroups holds so that   
the following theorem will be proved:
\proclaim{Theorem 1.1(Theorem 5.7)}
Let ${\frak L}_1$ and $ {\frak L}_2$ be 
$\lambda$-graph systems satisfying condition (I).
The following  are  equivalent:
 \roster
 \item There exists an isomorphism 
 $\Psi: \OLA \rightarrow  \OLB$
 such that $\Psi(\DLA) = \DLB$.
\item
$(X_{{\frak L}_1}, \sigma_{{\frak L}_1})$ 
and $(X_{{\frak L}_2}, \sigma_{{\frak L}_2})$ 
are continuously orbit equivalent.
\item 
There exists a homeomorphism 
$h: X_{{\frak L}_1} \rightarrow X_{{\frak L}_2}$ such that
$h \circ [\sigma_{{\frak L}_1} ]_{sc} \circ h^{-1} 
= [\sigma_{{\frak L}_2} ]_{sc}$. 
\endroster
\endproclaim
Let $\Lambda$ 
be the subshift presented by a $\lambda$-graph system $\frak L$
and $(X_\Lambda, \sigma_\Lambda)$ the right one-sided subshift for $\Lambda$.
There exists a natural factor map
$
\pi^{\frak L}_\Lambda: (X_{\frak L}, \sigma_{\frak L}) \longrightarrow 
(X_\Lambda,\sigma_\Lambda).
$
It induces an inclusion  
$C(X_\Lambda) \hookrightarrow C(X_{\frak L})$.
We regard the algebra $C(X_\Lambda)$ a subalgebra $\DLambda$
of $\DL$ and of $\OL$.
We say that two factor maps 
$\pi^{{\frak L}_1}_{\Lambda_1}$
and
$\pi^{{\frak L}_2}_{\Lambda_2}$
are 
continuously orbit equivalent
if there exist homeomorphisms
$h_{\frak L}:X_{{\frak L}_1}\longrightarrow X_{{\frak L}_1}$
and
$h_{\Lambda}:X_{\Lambda_1}\longrightarrow X_{\Lambda_2}$
such that
$\pi^{{\frak L}_2}_{\Lambda_2} \circ h_{\frak L} = h_\Lambda \circ \pi^{{\frak L}_1}_{\Lambda_1}$
and
there exist
continuous functions
$k_1,l_1:X_{{\frak L}_1}\longrightarrow \Zp$
and
$k_2,l_2:X_{{\frak L}_2}\longrightarrow \Zp$
such that 
$$
\align
\sigma_{{\frak L}_2}^{k_1(x)}(h_{\frak L} \circ \sigma_{{\frak L}_1}(x)) 
& = 
\sigma_{{\frak L}_2}^{l_1(x)}(h_{\frak L}(x)),
\quad
x \in X_{{\frak L}_1}, \\
\sigma_{{\frak L}_1}^{k_2(y)}(h_{\frak L}^{-1} \circ \sigma_{{\frak L}_2}(y)) 
& = 
\sigma_{{\frak L}_1}^{l_2(x)}(h_{\frak L}^{-1}(y)),
\quad
y \in X_{{\frak L}_2}.
\endalign
$$
Then we will prove
\proclaim{Theorem 1.2(Theorem 6.6)}
Let ${\frak L}_1$ and
${\frak L}_2$ be $\lambda$-graph systems satisfying condition (I)
and $\Lambda_1$ and $\Lambda_2$ their respect subshifts.
The following are equivalent:
\roster
\item"(i)"
There exists an isomorphism
$\Psi:\OLA \longrightarrow \OLB$ such that
$\Psi(\DLambdaA)=\DLambdaB$.
\item"(ii)"
The factor maps $\pi^{{\frak L}_1}_{\Lambda_1}$
and
$\pi^{{\frak L}_2}_{\Lambda_2}$
are
continuously orbit equivalent.
\item"(iii)" 
There exist homeomorphisms
$h_{\frak L}: X_{{\frak L}_1} \longrightarrow
 X_{{\frak L}_2}$
and
$h_\Lambda:X_{\Lambda_1} \longrightarrow X_{\Lambda_2}$
such that
$\pi^{{\frak L}_2}_{\Lambda_2} \circ h_{\frak L} = h_\Lambda \circ \pi^{{\frak L}_1}_{\Lambda_1}$
and
 $h_{\frak L}\circ [\sigma_{{\frak L}_1}]_{sc}\circ h_{\frak L}^{-1}
 = [\sigma_{{\frak L}_2}]_{sc}.
$
\endroster
\endproclaim

Let ${\frak L}^\Lambda$
be the canonical $\lambda$-graph system for $\Lambda$ (see [Ma2]).
Then the $C^*$-algebra $\OLambda$ coincides with the 
algebra
${\Cal O}_{{\frak L}^\Lambda}$.
The natural inclusion
$\iota: X_\Lambda \hookrightarrow X_{{\frak L}^\Lambda}$
induces a new topology on $X_\Lambda$.
The topological space is denoted by
$\widetilde{X}_\Lambda$.
Two subshifts
$(X_{\Lambda_1}, \sigma_{\Lambda_1})$ 
and 
$(X_{\Lambda_2},\sigma_{\Lambda_2})$
are said to be $\lambda$-{\it continuously orbit equivalent\/}
if
there exist a homeomorphism
$h:X_{\Lambda_1}\longrightarrow X_{\Lambda_2}$,
and
 continuous functions
$k_1,l_1:\widetilde{X}_{\Lambda_1}\longrightarrow \Zp$
and
$k_2,l_2:\widetilde{X}_{\Lambda_2}\longrightarrow \Zp$
such that 
$h$ is also homeomorphic from
$\widetilde{X}_{\Lambda_1}$ onto $\widetilde{X}_{\Lambda_2}$
such that
$$
\align
\sigma_{\Lambda_2}^{k_1(a)}(h\circ \sigma_{\Lambda_1}(a)) 
& =
\sigma_{\Lambda_2}^{l_1(a)}(h(a)), \qquad a \in X_{\Lambda_1}, \\
\sigma_{\Lambda_1}^{k_2(b)}(h^{-1}\circ \sigma_{\Lambda_2}(b)) 
& =
\sigma_{\Lambda_1}^{l_2(b)}(h^{-1}(b)), \qquad b \in X_{\Lambda_2}.
\endalign
$$
Then we will prove 
\proclaim{Theorem 1.3(Theorem 7.5)}
Let $\Lambda_1$ and $\Lambda_2$ be subshifts satisfying condition (I).
The following are equivalent:
\roster
\item
There exists an isomorphism
$\Psi:\OLambdaA \longrightarrow \OLambdaB$
such that
$\Psi(\DLambdaA) = \DLambdaB$.
\item
The subshifts 
$(X_{\Lambda_1},\sigma_{\Lambda_1})$
and
$(X_{\Lambda_2},\sigma_{\Lambda_2})$
are $\lambda$-continuously orbit equivalent.
\endroster
\endproclaim
The theorem is a generalization of a result in [Ma5] for topological Markov shifts. 

The results of this paper will be generalized to more general groupoid $C^*$-algebras in a forthcoming paper [Ma6]. 
Throughout the paper,
we denote by
$\Zp$ and $\Bbb N$ the set of nonnegative integers and the set of positive integers
respectively.
\medskip

\heading 2. Preliminaries
\endheading

Let $\frak L =(V,E,\lambda,\iota)$ be 
a $\lambda$-graph system 
 over $\Sigma$ with vertex set
$
V = \cup_{l \in \Zp} V_{l}
$
and  edge set
$
E = \cup_{l \in \Zp} E_{l,l+1}
$
that is labeled with symbols in $\Sigma$ by a map
$\lambda: E \rightarrow \Sigma$, 
and that is supplied with  surjective maps
$
\iota( = \iota_{l,l+1}):V_{l+1} \rightarrow V_l
$
for
$
l \in \Bbb \Zp.
$
Here the vertex sets $V_{l},l \in \Zp$
are finite disjoint sets.   
Also  
$E_{l,l+1},l \in \Zp$
are finite disjoint sets.
An edge $e$ in $E_{l,l+1}$ has its source vertex $s(e)$ in $V_{l}$ 
and its terminal  vertex $t(e)$ 
in
$V_{l+1}$
respectively.
Every vertex in $V$ has a successor and  every 
vertex in $V_l$ for $l\in \Bbb N$ has a predecessor. 
It is then required that there exists an edge in $E_{l,l+1}$
with label $\alpha$ and its terminal is  $v \in V_{l+1}$
 if and only if 
 there exists an edge in $E_{l-1,l}$
with label $\alpha$ and its terminal is $\iota(v) \in V_{l}.$
For 
$u \in V_{l-1}$ and
$v \in V_{l+1},$
put
$$
\align
E^{\iota}(u, v)
& = \{e \in E_{l,l+1} \ | \ t(e) = v, \iota(s(e)) = u \},\\
E_{\iota}(u, v)
& = \{e \in E_{l-1,l} \ | \ s(e) = u, t(e) = \iota(v) \}.
\endalign
$$ 
Then we require a bijective correspondence between 
$
E^{\iota}(u, v)
$
and
$
E_{\iota}(u, v)
$
that preserves labels
for each pair of vertices
$u, v$.
We call this property  the local property of $\frak L$. 
We henceforth assume that $\frak L$ is left-resolving, which means that 
$t(e)\ne t(f)$ whenever $\lambda(e) = \lambda(f)$ for $e,f \in E$.

Let $\Omega_{\frak L}$
be the compact Hausdorff space of the projective limit of the system
$\iota_{l,l+1}: V_{l+1}\rightarrow V_l, l \in \Zp,$ that is defined by
$$
\Omega_{\frak L} = \{ (v^l)_{l \in \Zp}\in 
\prod \Sb {l \in \Zp} \endSb V_l
 \ | \ \iota_{l,l+1}(v^{l+1}) = v^l, l\in \Zp \}.
$$
An element $v$ in $\Omega_{\frak L}$
is called  an $\iota$-orbit
or also a vertex.
Let $E_{\frak L}$ 
be the set of all triplets 
$(u, \alpha,v) \in \Omega_{\frak L} \times \Sigma \times \Omega_{\frak L}$,
where $u = (u^l)_{l \in \Zp}, v = (v^l)_{l \in \Zp}\in \Omega_{\frak L}$
such that
for each $l \in \Zp$, 
there exists $e_{l,l+1} \in E_{l,l+1}$ satisfying
$$
u^l = s(e_{l,l+1}),\quad
  v^{l+1} = t(e_{l,l+1}) \quad \text{ and } \quad
  \alpha = \lambda(e_{l,l+1}).
$$
Then the set
$
E_{\frak L}\subset \Omega_{\frak L} \times \Sigma \times \Omega_{\frak L}
$
is a  zero-dimensional continuous graph
 in the sense of Deaconu 
([Ma4;Proposition 2.1],[De],[De2],[De3],[De4]).
It has been also studied in [KM] as a Shannon graph.
 Following Deaconu [De2] and Krieger [Kr2], 
 we consider the set $X_{\frak L}$ 
 of all one-sided paths of $E_{\frak L}$:
$$
\align
X_{\frak L}= \{ (\alpha_n,u_n)_{n\in \Bbb N} \in 
\prod \Sb n\in \Bbb N \endSb
(\Sigma \times \Omega_{\frak L}) \ | \
& (u_{n},\alpha_{n+1},u_{n+1}) \in E_{\frak L} 
 \text{ for all } n\in \Bbb N \\ 
 \text{ and } 
& (u_{0},\alpha_1,u_{1}) \in E_{\frak L}
 \text{ for some } u_0 \in \Omega_{\frak L} \}.
\endalign
$$
The set
$X_{\frak L}$ becomes 
a zero-dimensional compact Hausdorff space
under the relative topology from the infinite product topology
of $\Sigma \times \Omega_{\frak L}$. 
For $x = (\alpha_n,u_n)_{n\in \Bbb N} \in
X_{\frak L}$,
the vertex $u_0 \in \Omega_{\frak L}$
satisfying
$(u_{0},\alpha_1,u_{1}) \in E_{\frak L}$
is unique
because $\frak L$ is left-resolving.
We denote it by $u_0(x)$. 
The shift map 
$
\sigma_{\frak L} :(\alpha_n,u_n)_{n\in \Bbb N} \in
X_{\frak L}
\rightarrow 
(\alpha_{n+1},u_{n+1})_{n\in \Bbb N}\in
X_{\frak L}
$
is 
a local homeomorphism 
by  [Ma4;Lemma 2.2].
We have a topological dynamical system 
$(X_{\frak L},\sigma_{\frak L})$
of a compact Hausdorff space $X_{\frak L}$ with a continuous surjection 
$\sigma_{\frak L}$.
The set
$$
X_\Lambda =
\{ (\alpha_n )_{n \in {\Bbb N}} \in \Sigma^{\Bbb N}
\mid   
(\alpha_n,u_n)_{n\in \Bbb N} \in
X_{\frak L} \}
$$
becomes the right one-sided subshift for the subshift 
$\Lambda$  
presented by $\frak L$
with shift transformation $\sigma_\Lambda$ defined by
$$
\sigma_\Lambda((\alpha_n )_{n \in {\Bbb N}}) 
= (\alpha_{n+1} )_{n \in {\Bbb N}},
\qquad
(\alpha_n )_{n \in {\Bbb N}} \in X_\Lambda.
$$
The factor map 
$$
\pi^{\frak L}_\Lambda : (\alpha_n,u_n)_{n\in \Bbb N} \in
X_{\frak L}
\rightarrow
(\alpha_n)_{n\in \Bbb N} \in
X_{\Lambda} 
$$
is a continuous surjective map satisfying
$$
\pi^{\frak L}_\Lambda \circ \sigma_{\frak L} = 
\sigma_\Lambda \circ \pi^{\frak L}_\Lambda.
$$
A word $\mu = \mu_1 \cdots \mu_k$ for $\mu_i \in \Sigma$
is said to be admissible for $X_\Lambda$ 
if $\mu$ appears in somewhere in some element $a$ in $X_\Lambda$.
 We denote by 
$B_k(X_\Lambda)$ 
the set of all admissible words of length $k \in \Bbb N$,
where $B_0(X_\Lambda)$ means the empty word $\emptyset$.
We set 
$B_*(X_\Lambda) = \cup_{k=0}^\infty B_k(X_\Lambda)$.
For $a = (a_n )_{n \in {\Bbb N}} \in X_\Lambda$ 
and positive integers $k,l$ with
$k\le l$, we put the word 
$a_{[k,l]} = (a_k,a_{k+1},\dots, a_l) \in B_{l-k +1}(X_\Lambda)$
and the right infinite sequence
$a_{[k,\infty)} =(a_k, a_{k+1}, \dots ) \in X_{\Lambda}$.
Similarly we use the notations
$B_k(X_{\frak L})$ defined by
the set 
$
\{(\alpha_n,u_n)_{n=1}^k \mid (\alpha_n,u_n)_{n\in \Bbb N} \in X_{\frak L}\}
$
and
$x_{[k,l]} = (x_k,\dots, x_l)$ 
for
 $x = (x_n )_{n \in {\Bbb N}} \in X_{\frak L}$.

Let us now briefly review the $C^*$-algebra
${\Cal O}_{\frak L}$
associated with $\lambda$-graph system $\frak L$.
The  $C^*$-algebras ${\Cal O}_{\frak L}$
 are generalization of the 
$C^*$-algebras associated with subshifts ([Ma4], cf.[CM]).
We denote by $\{v_1^l,\dots,v_{m(l)}^l\}$ 
the vertex set $V_l$.
Define the transition matrices $A_{l,l+1}, I_{l,l+1}$
of $\frak L$
by setting
for
$
i=1,2,\dots,m(l),\ j=1,2,\dots,m(l+1), \ \alpha \in \Sigma,
$ 
$$
\align
A_{l,l+1}(i,\alpha,j)
 & =
\cases
1 &  
    \text{ if } \ s(e) = v_i^l, \lambda(e) = \alpha,
                       t(e) = v_j^{l+1} 
    \text{ for some }    e \in E_{l,l+1}, \\
0           & \text{ otherwise,}
\endcases \\
I_{l,l+1}(i,j)
 & =
\cases
1 &  
    \text{ if } \ \iota_{l,l+1}(v_j^{l+1}) = v_i^l, \\
0           & \text{ otherwise.}
\endcases 
\endalign
$$
The $C^*$-algebra $\OL$
is realized as the universal unital $C^*$-algebra
generated by
partial isometries
$S_{\alpha}, \alpha \in \Sigma$
and projections
$E_i^l, i=1,2,\dots,m(l),\l\in \Zp 
$
 subject to the  following operator relations called $(\frak L)$:
$$
\align
\sum_{\alpha \in \Sigma} S_{\alpha}S_{\alpha}^*  & = 1, \tag 2.1 \\
 \sum_{i=1}^{m(l)} E_i^l   =  1, \qquad 
 E_i^l  &  =  \sum_{j=1}^{m(l+1)}I_{l,l+1}(i,j)E_j^{l+1}, \tag 2.2 \\
 S_\beta S_\beta^*E_i^l & =   E_i^{l} S_\beta S_\beta^*,
\tag 2.3 \\
S_{\beta}^*E_i^l S_{\beta}  =  
\sum_{j=1}^{m(l+1)}& A_{l,l+1}(i,\beta,j)E_j^{l+1},
\tag 2.4 
\endalign
$$
for $\beta \in \Sigma,$
$
i=1,2,\dots,m(l),\l\in \Zp. 
$
It is nuclear ([Ma4;Proposition 5.6]).
For a word $\mu = \mu_1\cdots\mu_k \in B_k(X_\Lambda)$,
we set 
$S_\mu = S_{\mu_1}\cdots S_{\mu_k}.$
The algebra of all finite linear combinations of the elements of the form
$$
S_{\mu}E_i^lS_{\nu}^* \quad \text{ for }\quad \mu,\nu \in B_*(X_\Lambda), \quad i=1,\dots,m(l),\quad l \in \Zp
$$
is a dense $*$-subalgebra of  $\OL$.
Let us denote by 
${\Cal A}_{\frak L}$ the $C^*$-subalgebra of $\OL$
generated by the projections
$E_i^l, i=1,\dots, m(l),\/ l \in \Zp$.
By the universality of the algebra $\OL$
 the algebra ${\Cal A}_{\frak L}$ 
is isomorphic to the commutative $C^*$-algebra $C(\Omega_{\frak L})$
of all complex valued continuous functions on $\Omega_{\frak L}$. 
We define  $C^*$-subalgebra
$\FKL $ with $k\le l$,
that is a finite dimensional algebra 
generated by 
$
          S_{\mu}E_i^lS_{\nu}^{*}, 
          \mu, \nu \in B_k(X_\Lambda), i=1,\dots,m(l).
$
Denote by 
$\FL$ the AF-subalgebra of $\OL$ generated by 
$\cup_{k,l}\FKL$.
For a vertex $v_i^l \in V_l$, 
put 
$$
\align
\Gamma^{+}(v_i^l) 
=\{ (\alpha_1,\alpha_2,\dots, ) \in \Sigma^{\Bbb N} & \ | \
 \text{ there exists an edge } e_{n,n+1} \in E_{n,n+1} 
 \text{ for } n \ge l  \\
 \text{such that }  v_i^l =  s(e_{l,l+1}),\ &
 t(e_{n,n+1})  = s(e_{n+1,n+2}), \ 
 \lambda(e_{n,n+1}) = \alpha_{n-l+1} \}
\endalign
$$
the set of all label sequences in $\frak L$ starting at $v_i^l$.
We say that $\frak L$ satisfies condition (I) if for each
$v_i^l \in V,$
the set $\Gamma^{+}(v_i^l)$ contains at least two distinct sequences.
Under the condition (I), 
the algebra $\OL$ can be realized as the  unique 
$C^*$-algebra subject to  the relations $(\frak L)$
([Ma4;Theorem 4.3]).
A $\lambda$-graph system $\frak L$ is said to be {\it irreducible \/}
if for a vertex $v \in V_l$  and an $\iota$-orbit
 $ x=(x_i)_{i\in \Zp} \in \Omega_{\frak L},$
 there exists a $\lambda$-path starting at $v$ and terminating at
 $x_{l+N}$ for some $N\in \Bbb N$.
If $\frak L$ is irreducible with  condition (I),
the $C^*$-algebra $\OL$ is simple
([Ma4;Theorem 4.7]).

Let
$
\DL
$
be the $C^*$-subalgebra of $\FL$ generated by 
$S_\mu E_i^l S_\mu^*, \enskip
\mu \in B_*(X_\Lambda), i=1,\dots,m(l), l \in \Zp
$
and
$
\DLambda 
$
the $C^*$-subalgebra of $\DL$
generated by
$S_\mu S_\mu^*, \enskip \mu \in B_*(X_\Lambda).$ 
For 
$\mu = \mu_1\cdots\mu_k \in B_k(X_\Lambda)$
and 
$v_i^l \in V_l$,
we set the cylinder set
$$
U_{\mu,v_i^l} =
\{ (\alpha_n,u_n) \in X_{\frak L} 
\mid \alpha_1 = \mu_1,\dots,\alpha_1 =\mu_k, \/ u_k^l = v_i^l \}
$$
of $X_{\frak L}$
where
$u_k = (u_k^l)_{l \in \Zp} \in \Omega_{\frak L}$.  
Let
$\chi_{U_{\mu,v_i^l}}$ denote the chracteristic function
on
$X_{\frak L}$ for the cylinder set 
$U_{\mu,v_i^l}$.
Then the correspondence
$ 
S_\mu E_i^l S_\mu^* \in \DL
 \longleftrightarrow \chi_{U_{\mu,v_i^l}} \in C(X_{\frak L})
$
yields an isomorphism between 
$\DL$ and
$C(X_{\frak L}).
$
Similarly 
let
$
U_\mu = \{ (a_n )_{n \in \Bbb N} \in X_{\Lambda} 
\mid a_1 =\mu_1,\dots, a_k = \mu_k \}
$
be the cylinder set of $X_\Lambda$.
The correspondence
$ 
S_\mu S_\mu^* \in \DLambda \longleftrightarrow \chi_\mu \in C(X_{\Lambda})
$
yields an isomorphism between 
$ \DLambda$
and
$C(X_\Lambda)$.

By the universality for the relations $(\frak L)$, 
the correspondence
$
S_\alpha \longrightarrow e^{\sqrt{-1}t}S_\alpha, \alpha \in \Sigma, \,
E_i^l \longrightarrow E_i^l, i=1,\dots,m(l), \, l \in \Zp$ 
for
$
e^{\sqrt{-1}t} \in 
{\Bbb T}=\{ e^{\sqrt{-1}t} \mid t \in [0,2 \pi]\}
$
gives rise to  an action 
$
\rho: {\Bbb T} \rightarrow \Aut(\OL)
$
called gauge action.
The fixed point algebra of $\OL$ under $\rho$
is the AF-algebra $\FL$.
We denote by $E: \OL \rightarrow \FL$
the conditional expectation defined by
$E(a) = \int_{\Bbb T}\rho_t(a) dt $ for $a\in \OL$.

The following lemma is basic in our further discussions.
\proclaim{Lemma 2.1( [Ma3;Proposition 3.3], cf.[CK;Remark 2.18])}
Suppose that $\frak L$ satisfies condition (I). 
Then we have
$\DLambda^{\prime} \cap \OL = \DL$
and hence $\DL^{\prime} \cap \OL = \DL.$
\endproclaim
This means that 
the algebra $\DL$ is maximal abelian in $\OL$.
\demo{Proof}
The proof of $\DLambda^{\prime} \cap \OL = \DL$
is completely similar to the proof of [Ma3;Proposition 3.3].
Since
$\DL \subset \DL^{\prime}\cap \OL \subset \DLambda\cap \OL$,
we have $\DL^{\prime} \cap \OL = \DL.$
\qed
\enddemo
In [Ma5], 
a representation of the Cuntz-Krieger algebra 
$\OA$ on a Hilbert space having the shift space $X_A$ 
as a complete orthonormal basis has been used. 
Let us generalize the representation to the $C^*$-algebras $\OL$ as in the following way.
 Let
${\frak H}_{\frak L}$  
be the Hilbert space with its complete orthonormal system
$ e_x , x \in X_{\frak L}$. 
The Hilbert space is not separable.
Consider the partial isometries
$T_\alpha:
{\frak H}_{\frak L} \rightarrow {\frak H}_{\frak L}, \alpha \in \Sigma
$ 
and projections
$P_i^l: 
{\frak H}_{\frak L} \rightarrow {\frak H}_{\frak L}, i=1,\dots,m(l)
$ 
defined by
$$
T_\alpha e_x = 
\cases
e_{y} & 
\text{ if there exists an }\iota\text{-orbit } u_{-1}\in \Omega_{\frak L} ;
(u_{-1},\alpha,u_0(x)) \in E_{\frak L},\\
0 & \text{ otherwise }  
\endcases
$$
where $y = ( (\alpha,u_0(x)), (\alpha_1,u_1), (\alpha_2,u_2),\dots ) 
\in X_{\frak L}$ for
$x = ( (\alpha_1,u_1), (\alpha_2,u_2),\dots ) 
\in X_{\frak L}$
and
$$
P_i^l e_x =
\cases
e_x & \text{ if } u_0(x)^l = v_i^l, \\
0 & \text{ otherwise }  
\endcases
$$
where
$u_0(x) = (u_0(x)^l)_{l \in \Zp}\in \Omega_{\frak L}.$
\proclaim{Lemma 2.2}
The partial isometries
$T_\alpha, \alpha \in \Sigma
$ 
and the projections
$P_i^l, i=1,\dots,m(l)
$ 
on the Hilbert space
$
{\frak H}_{\frak L}
$
satisfy the relation $(\frak L)$.
Hence if $\frak L$ satisfies condition (I), 
the correspondence
$S_\alpha \rightarrow T_\alpha$ and
$E_i^l\rightarrow P_i^l$
gives rise to a faithful representation of the $C^*$-algebra
$\OL$ on ${\frak H}_{\frak L}$.
\endproclaim
We call it
the universal shift representation 
of $\OL$ on ${\frak H}_{\frak L}$.
In what follows, we assume that $\frak L$ satisfies condition (I)
and regard the algebra $\OL$ as the $C^*$-algebra generated by 
$T_\alpha, \alpha \in \Sigma
$ 
and 
$P_i^l, i=1,\dots,m(l)
$ 
on the Hilbert space
$
{\frak H}_{\frak L}.
$
 
\heading 3. Topological full inverse semigroups
\endheading
For $x = (x_n )_{n \in \Bbb N} \in X_{\frak L}$,
the orbit $orb_{\sigma_{\frak L}}(x)$ of $x$ is defined by
$$
orb_{\sigma_{\frak L}}(x) 
= 
\cup_{k=0}^\infty \cup_{l=0}^\infty 
\sigma_{\frak L}^{-k}(\sigma_{\frak L}^l(x)) \subset X_{\frak L}.
$$
Hence  
$ y =( y_n )_{n \in \Bbb N} \in X_{\frak L}$ 
belongs to $orb_{\sigma_{\frak L}}(x)$ 
if and only if
there exists a a finite sequence 
$z_1 \cdots z_k \in B_k(X_{\frak L})$
such that 
$$
y = (z_1,\dots, z_k, x_{l+1}, x_{l+2},\dots )
\qquad \text{for some } k, l \in \Zp.
$$
We denote by $\Homeo(X_{\frak L})$
the group of all homeomorphisms on $X_{\frak L}$.
We define the full group 
$[\sigma_{\frak L}]$ and the topological full group
$[\sigma_{\frak L}]_c$
for $(X_{\frak L},\sigma_{\frak L})$
as in the following way.

\noindent
{\bf Definition.}
Let
$
[\sigma_{\frak L} ] 
$ be the set of all homeomorphism
$\tau \in \Homeo(X_{\frak L})
$ 
such that
$
\tau(x) \in  orb_{\sigma_{\frak L}}(x)
$ 
for all
$ 
x \in X_{\frak L}.
$ 
We call $[\sigma_{\frak L}]$ 
the full group of $(X_{\frak L},\sigma_{\frak L})$.
Let $[\sigma_{\frak L} ]_c$ be the set of all $\tau$ in $[\sigma_{\frak L}]$ such that  
there exist continuous functions 
$k, l : X_{\frak L} \rightarrow \Zp$ 
such that 
$$
\sigma_{\frak L}^{k(x)}(\tau(x) )=\sigma_{\frak L}^{l(x)}(x)
\quad\text{ 
for all } x \in X_{\frak L}. \tag 3.1
$$
We call   
$[\sigma_{\frak L} ]_c$ 
the topological full group for $(X_{\frak L},\sigma_{\frak L})$.

If a subshift is not a sofic shift,
 the full groups 
are not necessarily large enough to cover the orbit structure.
Hence to study of orbit structure of general subshifts,
we will extend the notion of full groups to full inverse semigroups
as in the following way.
Let $\tau:U \rightarrow V$ be a homeomorphism
from a clopen set $U \subset X_{\frak L}$ 
onto a clopen set $V \subset X_{\frak L}$.
We call $\tau$ 
a partial homeomorphism. 
Let us denote by $X_{\tau}$ and $Y_{\tau}$
the clopen sets $U$ and $V$ respectively.
We denote by 
$PH(X_{\frak L})$ the set of all partial homeomorphisms of $X_{\frak L}$.
Then $PH(X_{\frak L})$ has a natural structure of inverse semigroup
(cf. [Pat]).
We define the full inverse semigroup 
$[\sigma_{\frak L}]_s$ 
and the topological full inverse semigroup
$[\sigma_{\frak L}]_{sc}$
for $(X_{\frak L},\sigma_{\frak L})$
as in the following way.

\noindent
{\bf Definition.}
Let
$
[\sigma_{\frak L} ]_{s} 
$ be the set of all partial homeomorphisms
$\tau \in PH(X_{\frak L})
$ 
such that
$
\tau(x) \in  orb_{\sigma_{\frak L}}(x)
$ 
for all
$ 
x \in X_{\tau}.
$ 
We call $[\sigma_{\frak L}]_s$ 
the full inverse semigroup of $(X_{\frak L},\sigma_{\frak L})$.
Let $[\sigma_{\frak L} ]_{sc}$ 
be the set of all $\tau$ in $[\sigma_{\frak L}]_s$ 
such that  
there exist continuous functions 
$k, l : X_{\tau} \rightarrow \Zp$ 
such that 
$$
\sigma_{\frak L}^{k(x)}(\tau(x) )=\sigma_{\frak L}^{l(x)}(x)
\quad\text{ 
for all } x \in X_{\tau}. \tag 3.2
$$
We call   
$[\sigma_{\frak L} ]_{sc}$ 
the topological full inverse semigroup for 
$(X_{\frak L},\sigma_{\frak L})$.
The maps $k,l$ above are called orbit cocycles for $\tau$,
and sometimes written as $k_\tau, l_\tau$ respectively.
We remark that the orbit cocyles are not necessarily uniquely determined  
for $\tau$.
It is clear that 
$[\sigma_{\frak L}]_s$ is a subsemigroup of $PH(X_{\frak L})$ 
and 
$[\sigma_{\frak L}]_{sc}$ is a subsemigroup of $[\sigma_{\frak L}]_c$. 
Although $\sigma_{\frak L}$ does not belong to $[\sigma_{\frak L}]_{sc}$,
the following lemma shows that $\sigma_{\frak L}$ locally belongs to
$[\sigma_{\frak L}]_{sc}$, and  
that
 $[\sigma_{\frak L}]_{sc}$ 
is not trivial in any case.
\proclaim{Lemma 3.1}
For any $\mu= (\mu_1,\dots, \mu_k)  \in B_k(X_{\Lambda})$
and
$v_i^l \in V_l$ with $2 \le k \le l$,
there exists $\tau_{\mu, v_i^l} \in [\sigma_{\frak L}]_{sc}$
such that  
$$
 \tau_{\mu, v_i^l}(x) = 
\sigma_{\frak L}(x) \qquad
 \text{ for } x \in U_{\mu,v_i^l}.
\tag 3.3
$$
\endproclaim
\demo{Proof}
Put  
$\nu = (\mu_2,\dots, \mu_k) \in B_{k-1}(X_{\Lambda})$.
Then the map
$
 \tau_{\mu, v_i^l}: U_{\mu,v_i^l} \longrightarrow   U_{\nu,v_i^l}
$
defined by
$ \tau_{\mu, v_i^l}(x) = 
\sigma_{\frak L}(x)
$ 
for
$ x \in U_{\mu,v_i^l}$
is a partial homeomorphism,
and it belongs to $[\sigma_{\frak L}]_{sc}$.
\qed
\enddemo
\proclaim{Lemma 3.2}
For 
$
x = (x_n)_{n \in \Bbb N}\in X_{\frak L}
$ 
with 
$ 
x_n = (\alpha_n,u_n), n\in {\Bbb N},
$ 
put $u_0 = u_0(x) \in \Omega_{\frak L}$.
Let $\alpha_0 \in \Sigma$ be a symbol such that 
$(\alpha_{n-1},u_{n-1})_{n \in \Bbb N} \in X_{\frak L}$.
Then there exists 
$\tau \in [\sigma_{\frak L}]_{sc}$
with a clopen set $X_{\tau}\subset X_{\frak L}$
such that
$ x \in X_{\tau}$ and
$\tau(y) = (y_{n-1})_{n\in {\Bbb N}}$ 
for all
$ y= (y_n)_{n\in {\Bbb N}}\in X_\tau$,
where
$y_0 = (\alpha_0,u_0(y))$.
\endproclaim
\demo{Proof}
Let
$
X_\tau
$ 
be the clopen set 
$
 U_{\mu,v_i^l} 
$ for
$\mu = \alpha_1 \alpha_2 \in B_2(X_\Lambda)$
and
$v_i^l = u_2^2 \in V_2$,
where
$
u_2 =(u_2^l)_{l\in \Zp}\in \Omega_{\frak L},
$
so that 
$x$ belongs to $X_\tau $.
One has 
$(y_{n-1})_{n \in {\Bbb N}} \in X_{\frak L}$
for
$(y_n)_{n \in {\Bbb N}} \in X_\tau$,
where
$y_0 = (\alpha_0,u_0(y))$.
By setting $\tau(y) = (y_{n-1})_{n\in \Bbb N}$
for
$ y = (y_n)_{n \in \Bbb N} \in X_{\frak L}$,
we have
$\sigma_{\frak L}(\tau(y)) = y$ for $ y \in X_\tau$
so that 
$\tau \in [\sigma_{\frak L}]_{sc}$.  
\qed
\enddemo
For $x \in X_{\frak L}$, put
$
[\sigma_{\frak L}]_{sc} (x) 
= \{ \tau(x) \in X_{\frak L} \mid \tau \in [\sigma_{\frak L}]_{sc} 
\text{ with }   X_\tau \ni x \}.
$
\proclaim{Lemma 3.3}
$
[\sigma_{\frak L}]_{sc} (x) = orb_{\sigma_{\frak L}}(x).
$
\endproclaim
\demo{Proof}
For any $\tau \in [\sigma_{\frak L}]_{sc}$
with
$  X_\tau \ni x$, 
one sees 
$\tau(x) \in orb_{\sigma_{\frak L}}(x)$
and hence
$
[\sigma_{\frak L}]_{sc} (x) \subset orb_{\sigma_{\frak L}}(x).
$
For the other inclusion relation,
by the previous lemmas,
for $x = (x_n )_{n \in \Bbb N} \in X_{\frak L}$
and $x_0 =(\alpha_0,u_0(x)) \in \Sigma \times \Omega_{\frak L}$,
there exist
$\tau_1, \tau_2 \in [\sigma_{\frak L}]_{sc}$
such that
$$
\tau_1(x) = (x_{n-1})_{n \in \Bbb N}, \qquad
\tau_2(x) = (x_{n+1})_{n \in \Bbb N} \in X_{\frak L}
$$
so that both
$
 (x_{n-1})_{n \in \Bbb N}
$ and
$
 (x_{n+1})_{n \in \Bbb N} 
$ belong to
$[\sigma_{\frak L}]_{sc}(x)$.  
Since $[\sigma_{\frak L}]_c$ is a semigroup, 
one sees that
$$
[\sigma_{\frak L}]_{sc}(x) \ni (x_{-k}, \dots, x_{-1},
x_0, x_{l+1},x_{l+2}, \dots, )
$$
for all $k,l\in \Zp$
with
$
(x_{-k}, \dots, x_{-1},x_0, x_{l+1},x_{l+2}, \dots ) \in X_{\frak L}.
$
Hence
$
[\sigma_{\frak L}]_{sc} (x) \supset orb_{\sigma_{\frak L}}(x).
$
\qed
\enddemo

\heading 4. Full inverse semigroups and normalizers  
\endheading
Let us denote by
${\Cal U}(\OL)$ 
the group of unitaries of $\OL$ 
and
${\Cal U}(\DL)$ 
the group of unitaries of
$\DL$ respectively.
As in [Ma5],
the topological full group $[\sigma_{\frak L}]_{c}$
will correspond to 
the normalizer $N(\OL,\DL)$ of $\DL$ in $\OL$
defined by
$$
N(\OL,\DL) = \{ v \in {\Cal U}(\OL) \mid v \DL v^* = \DL \}.
$$
For   
the topological full inverse semigroup $[\sigma_{\frak L}]_{sc}$,
we will define 
the normalizer 
$N_s(\OL,\DL)$ of partial isometries
as in the following way:
$$
N_s(\OL,\DL) = \{ v \in \OL 
\mid v \text{ is a partial isometry };
v \DL v^* \subset \DL, v^* \DL v \subset \DL \}.
$$
 It is easy to see that 
$
N_s(\OL,\DL)
$ 
has a natural structure of inverse semigroup.
We will identify  the subalgebra 
$\DL$ of $\OL$ with the algebra $C(X_{\frak L})$.  
For a partial isometry $v \in \OL$,
put $Ad(v)(x) = v x v^*$ for $x \in \OL$. 
The following proposition holds. 
\proclaim{Proposition 4.1}
For $\tau \in [\sigma_{\frak L}]_{sc}$, 
there exists a partial isometry 
$u_\tau \in N_s(\OL,\DL)$
such that
$$
Ad(u_\tau)(f) = f \circ \tau^{-1}
\quad \text{ for } f \in C(X_\tau),
\qquad 
Ad(u_\tau^*)(g) = g \circ \tau
\quad
\text{ for } g \in C(Y_\tau),
$$
and the correspondence 
$\tau \in [\sigma_{\frak L}]_{sc} \longrightarrow u_\tau \in N_s(\OL,\DL)$  
is a homomorphism of inverse semigroup. 
If in particular
$\tau \in [\sigma_{\frak L}]_c$,
the partial isometry $u_\tau$ is a unitary so that 
$u_\tau \in N(\OL,\DL)$.
\endproclaim
\demo{Proof}
Let the $C^*$-algebra $\OL$ be represented on the Hilbert space
${\frak H}_{\frak L}$ with complete orthonormal basis 
$\{ e_x\mid x \in X_{\frak L}\}$.
Put the subspaces
$$
{\frak H}_{X_\tau} = span\{ e_x \mid x \in X_\tau \},
\qquad  
{\frak H}_{Y_\tau} = span\{ e_x \mid x \in Y_\tau \}.
$$
Since
$\tau:X_\tau \longrightarrow Y_\tau$
is a homeomorphism,
the operator  
$u_\tau:  
{\frak H}_{X_\tau}\longrightarrow
{\frak H}_{Y_\tau}
$
defined by
$u_\tau{(e_x) }= e_{\tau(x)}$
for $x \in X_\tau$
yields a partial isometry on ${\frak H}_{\frak L}$.
By a similar manner to the proof of 
[Ma5:Proposition 4.1],
one knows that 
$u_\tau$ belongs to
$N_s(\OL,\DL)$. 
\qed.
\enddemo
For $v \in N_s(\OL,\DL)$,
put the projections
$p_v = v^* v, q_v = v v^*$ 
in $\DL$,
and
the clopen subsets
Let $X_v = \supp(p_v), Y_v= \supp(q_v)$
of 
$X_{\frak L}$.
Then
$Ad(v):\DL p_v \longrightarrow \DL q_v$ is an isomorphism 
and induces a partial homeomorphism
$\tau_v: X_v \longrightarrow Y_v$
such that 
$$
Ad(v)(f) = f \circ \tau_v^{-1}
\quad \text{ for } f \in C(X_v),
\qquad 
Ad(v^*)(g) = g \circ \tau_v
\quad
\text{ for } g \in C(Y_v).
$$
We will prove that $\tau_v$ 
gives rise to an element of $[\sigma_{\frak L}]_{sc}$.
Since the proof basically follows a line of the proof of 
[Ma5:Proposition 4.7], 
we will give a sketch of the proof.
Fix $v\in N_s(\OL,\DL)$ for a while.
\proclaim{Lemma 4.2}
\roster
\item"(i)"
There exists a family
$v_m, m\in {\Bbb Z}$ of partial isometries in $\OL$
such that all but finitely many $v_m, m\in {\Bbb Z}$ are zero, 
and
{\roster
\item
$ v = \sum_{m \in {\Bbb Z}} v_m :$ finite sum. 
\item $v_m^*v_m, v_mv_m^*$ are projections in $\DL$ for $m \in \Bbb Z$. 
\item $v_m \DL v_m^* \subset \DL$ and $v_m^* \DL v_m \subset \DL$ 
      for $m \in \Bbb Z$.
\item $v_m^* v_{m'} = v_m v_{m'}^* = 0$ for $m \ne {m'}$.
\item $v_0 \in \FL$.
\endroster}
\item"(ii)"
For a fixed $n \in {\Bbb N}$, 
there exist partial isometries
$ v_\mu , v_{-\mu} \in \FL$ for each
 $\mu \in B_n(X_\Lambda)$
satisfying the following conditions:
{\roster
\item
$ v_n = \sum_{\mu \in B_n(X_\Lambda)} S_\mu v_{\mu}$
and
$ v_{-n} = \sum_{\mu \in B_n(X_\Lambda)}  v_{-\mu}S_\mu^*$.
\item
 $v_{\mu}^*v_\mu$, $S_\mu v_\mu v_\mu^*S_\mu^*$, 
 $ S_\mu v_{-\mu}^*v_{-\mu}S_\mu^* $ 
and $v_{-\mu} v_{-\mu}^*$ are projections in $\DL$ such that
$$
\align
v_n^* v_n & = \sum_{\mu \in B_n(X_\Lambda)} v_\mu^* v_\mu, \qquad
v_n v_n^* = \sum_{\mu \in B_n(X_\Lambda)} S_\mu v_\mu v_\mu^* S_\mu^*,\\
v_{-n}^* v_{-n} &
 = \sum_{\mu \in B_n(X_\Lambda)} S_\mu v_{-\mu}^* v_{-\mu}S_\mu^*, \qquad
v_{-n} v_{-n}^* = \sum_{\mu \in B_n(X_\Lambda)} v_{-\mu} v_{-\mu}^*.
\endalign
$$
\item
$v_{\mu}v_\nu^* =  v_{-\mu}^* v_{-\nu}=0$ 
for $ \mu, \nu \in B_n(X_\Lambda)$ with 
$\mu \ne \nu$.
\item
The algebras 
$v_\mu \DL v_\mu^*, v_\mu^* \DL v_\mu, v_{-\mu} \DL v_{-\mu}^*$ 
and $v_{-\mu}^* \DL v_{-\mu}$ are contained in $\DL$. 
\endroster}
\endroster
\endproclaim
\demo{Proof}
(i)
Put a partial isometry $g(t) = v^* \rho_t(v) \in \OL$ for $t \in {\Bbb T}$.
For $f \in \DL$, it follows that          
$
\rho_t(v) f \rho_t(v)^* = \rho_t(v f v^* ) = v f v^*
$
and hence
$$
g(t) f =  v^*\rho_t(v) f \rho_t(v^*)\rho_t(v)=v^*vfv^*\rho_t(v)=fg(t)
$$ 
so that 
$g(t)$
 commutes with each element of $\DL$.
By Lemma 2.1, $g(t)$ belongs to the algebra $\DL$.
Since 
$g(t)^* = g(-t)$ and $g(t+s) = g(t) g(s)$,
by putting
$$
v_m = \int_{\Bbb T}\rho_t(v) e^{-\sqrt{-1}mt} dt, 
\qquad 
\hat{g}(m) = \int_{\Bbb T} g(t) e^{-\sqrt{-1}mt} dt
\quad \text{ for } m \in {\Bbb Z}.
$$
one has $ v_m = v \hat{g}(m)$.  
By a similar argument to the proof of [Ma5:Lemma 4.2],
one has the assertions (1),(2),(3), (4)  and (5).

(ii)
Put for $ \mu \in B_n(X_{\frak L})$,
$$
v_\mu = E(S_\mu^* v), \qquad v_{-\mu} =E(v S_\mu). 
$$
By a similar argument to the proof of [Ma5:Lemma 4.3],
one has the assertions (1),(2),(3) and (4).
\qed
\enddemo
For $u \in N_s(\OL,\DL)$,
let
$\tau_u: X_u \rightarrow Y_u$ 
be the induced homeomorphism.
\proclaim{Lemma 4.3}
Keep the above notation.
For $x = (x_n)_{n \in {\Bbb N}}\in X_u$
with
$x_n = (\alpha_n,u_n(x))$,
$ u_n(x) =(u_n^l(x))_{l\in \Zp}$,
put $ y =(y_n)_{n \in {\Bbb N}} =\tau_u(x) \in Y_u$, 
where
$y_n = (\beta_n,u_n(y))$,
$ u_n(y) =(u_n^l(y))_{l\in \Zp}$.
For a fixed integer $l\in \Zp$, 
take
$i(x_n) \in \{1,\dots,m(l)\}$
and
$i(y_n) \in \{1,\dots,m(l)\}$
such that 
$v_{i(x_n)}^l =u_n^l(x)$
and
$v_{i(y_n)}^l =u_n^l(y)$
respectively.
Then we have
$$
\| E_{i(y_n)}^l S^*_{\beta_1 \cdots \beta_n} u S_{\alpha_1 \cdots \alpha_n}
 E_{i(x_n)}^l\| = 1
\quad \text{ for all } n \in {\Bbb N}. 
$$
\endproclaim
\demo{Proof}
It suffices to show that
$
E_{i(y_n)}^l S^*_{\beta_1 \cdots \beta_n} u S_{\alpha_1 \cdots \alpha_n}
 E_{i(x_n)}^l \ne 0.
$
Since $v_{i(y_n)}^l = u_n^l(y)$,
one sees that 
$E_{i(y_n)}^l e_{{\sigma_{\frak L}}^n(y)} = e_{{\sigma_{\frak L}}^n(y)}$
so that
$$
\align
&(E_{i(y_n)}^l S^*_{\beta_1 \cdots \beta_n} u S_{\alpha_1 \cdots \alpha_n}
 E_{i(x_n)}^lS_{\alpha_1 \cdots \alpha_n}^*u^*S_{\beta_1 \cdots \beta_n}
E_{i(y_n)}^l 
e_{{\sigma_{\frak L}}^n(y)} \mid  e_{{\sigma_{\frak L}}^n(y)})\\
=& ( Ad(u)(S_{\alpha_1 \cdots \alpha_n}
 E_{i(x_n)}^lS_{\alpha_1 \cdots \alpha_n}^*)S_{\beta_1 \cdots \beta_n}
e_{{\sigma_{\frak L}}^n(y)} \mid S_{\beta_1 \cdots \beta_n}
e_{{\sigma_{\frak L}}^n(y)})\\ 
 =& ( Ad(u)(S_{\alpha_1 \cdots \alpha_n}
 E_{i(x_n)}^lS_{\alpha_1 \cdots \alpha_n}^*)
e_{y} \mid e_{y}).
\endalign 
$$
Consider the cylinder set
$$
U_{\alpha_1\cdots\alpha_n,v_{i(x_n)}^l}
=\{ (\gamma_m,u_m)_{m \in \Bbb N} \in X_{\frak L}
\mid \gamma_1 = \alpha_1, \dots, \gamma_n = \alpha_n, u_n^l = v_{i(x_n)}^l \}
$$
of $X_{\frak L}$.
As
$S_{\alpha_1 \cdots \alpha_n}
 E_{i(x_n)}^lS_{\alpha_1 \cdots \alpha_n}^* = 
\chi_{U_{\alpha_1\cdots\alpha_n,v_{i(x_n)}^l}}
$
and
$$
Ad(u)(\chi_{U_{\alpha_1\cdots\alpha_n,v_{i(x_n)}^l}})e_y  
= (\chi_{U_{\alpha_1\cdots\alpha_n,v_{i(x_n)}^l}}\circ \tau_u^{-1})(y) e_y
= \chi_{U_{\alpha_1\cdots\alpha_n,v_{i(x_n)}^l}}(x) e_y
=e_y,
$$
we have
$$
(E_{i(y_n)}^l S^*_{\beta_1 \cdots \beta_n} u S_{\alpha_1 \cdots \alpha_n}
 E_{i(x_n)}^l S_{\alpha_1 \cdots \alpha_n}^*u^*S_{\beta_1 \cdots \beta_n}
E_{i(y_n)}^l 
e_{{\sigma_{\frak L}}^n(y)} \mid  e_{{\sigma_{\frak L}}^n(y)})
= (e_y \mid e_y ) =1
$$
so that
$
E_{i(y_n)}^l S^*_{\beta_1 \cdots \beta_n} u S_{\alpha_1 \cdots \alpha_n}
 E_{i(x_n)}^l \ne 0.
$
\qed
\enddemo
\proclaim{Lemma 4.4} 
Keep the above situation.
Assume in particular that $u \in \FL$.
Then  there exists $k \in \Bbb N$ such that
for  all 
$
x=(x_n )_{n \in {\Bbb N}} \in X_u
$
$$
\tau_u(x)_n = x_n
\qquad
\text{ for all }n > k
$$
where $\tau_u(x) = (\tau_u(x)_n)_{n \in \Bbb N}$.
\endproclaim
\demo{Proof} 
Suppose that for any $k \in {\Bbb N}$
there exist  
$x \in X_u$ and $N >k$ such that
$\tau_u(x)_N \ne x_N$.
Put $y_n = \tau_u(x)_n, n \in \Bbb N$.
Now $u \in \FL$ 
so that take $u_{0} \in {\Cal F}_{l_0}^{k_0}$
 for some $k_0\le l_0$
such that
$ \| u - u_{0} \| < \frac{1}{2}.$
Take $ x \in X_u$ and $N_0 > k_0$ 
such as $y_{N_0} \ne x_{N_0}$.
Since
$
x_{N_0} = (\alpha_{N_0}, u_{N_0}(x)),
y_{N_0} = (\beta_{N_0}, u_{N_0}(y))
$
and
$u_{N_0}(x) = (u_{N_0}^l(x))_{l\in \Bbb N},
u_{N_0}(y) = (u_{N_0}^l(y))_{l\in \Bbb N}\in \Omega_{\frak L},
$
one has
$\alpha_{N_0} \ne \beta_{N_0}$
or
there exists $l_1$ 
such that 
$ 
u_{N_0}^l(x)
\ne u_{N_0}^l(y)
$
fo all $l \ge l_1$.
As 
$ 
u_{N_0}^l(x) = v_{i(x_{N_0})}^l,
u_{N_0}^l(y) = v_{i(y_{N_0})}^l,
$
the later condition is equivalent to the condition that
$ 
E_{i(x_{N_0})}^l
\ne E_{i(y_{N_0})}^l
$
fo all $l \ge l_1$.
Now
 $u_{0}\in {\Cal F}_{l_0}^{k_0} \subset {\Cal F}_{l_0'}^{N_0-1}$,
where
$l_0' = l_0 + N_0 -1-k_0$,
it is written 
as 
$$
u_{0} = \sum_{\xi, \eta \in B_{N_0-1}(X_{\Lambda}), j=1,\dots,m(l_0')} 
c_{\xi,j,\eta}S_\xi E_j^{l_0'}S_\eta^* \in {\Cal F}_{l_0'}^{N_0-1}
\quad \text{ for some } c_{\xi,j, \eta} \in {\Bbb C}.
$$
Hence we have
$$
\align
& S_{\beta_1\cdots\beta_{N_0-1}}^* u_0 S_{\alpha_1\cdots\alpha_{N_0-1}}\\
=&
\sum_{j=1}^{m(l_0')}
c_{\beta_1\cdots\beta_{N_0-1},j, \alpha_1\cdots\alpha_{N_0-1}} 
S_{\beta_1\cdots\beta_{N_0-1}}^* S_{\beta_1\cdots\beta_{N_0-1}}E_j^{l_0'}S_{\alpha_1\cdots\alpha_{N_0-1}}^*
S_{\alpha_1\cdots\alpha_{N_0-1}}.
\endalign
$$
Take an integer 
$l_1' $ such that 
$l_1' \ge \max\{l_1,l_0'\}$
and hence the condition
$\alpha_{N_0} \ne \beta_{N_0}$
or
$E_{i(x_{N_0})}^{l_1'} \cdot E_{i(y_{N_0})}^{l_1'} =0$
holds.
It follows that
$$
\align
 &
E_{i(y_{N_0})}^{l_1'} S^*_{\beta_1 \cdots \beta_{N_0}}
 u_0 S_{\alpha_1 \cdots \alpha_{N_0}}
E_{i(x_{N_0})}^{l_1'} \\
=& 
\sum_{j=1}^{m(l_0')}
c_{\beta_1\cdots\beta_{N_0-1},j, \alpha_1\cdots\alpha_{N_0-1}} 
E_{i(y_{N_0})}^{l_1'} 
S_{\beta_1\cdots\beta_{N_0}}^* S_{\beta_1\cdots\beta_{N_0-1}}E_j^{l_0'}S_{\alpha_1\cdots\alpha_{N_0-1}}^*
S_{\alpha_1\cdots\alpha_{N_0}}
E_{i(x_{N_0})}^{l_1'}. 
\endalign
$$
Since 
$
S_{\beta_1\cdots\beta_{N_0-1}}^* S_{\beta_1\cdots\beta_{N_0-1}}E_j^{l_0'}S_{\alpha_1\cdots\alpha_{N_0-1}}^*
S_{\alpha_1\cdots\alpha_{N_0-1}} 
$ 
belongs to
$\DL$,
one has
$$
E_{i(y_{N_0})}^{l_1'} 
S_{\beta_1\cdots\beta_{N_0}}^* S_{\beta_1\cdots\beta_{N_0-1}}E_j^{l_0'}S_{\alpha_1\cdots\alpha_{N_0-1}}^*
S_{\alpha_1\cdots\alpha_{N_0}}
E_{i(x_{N_0})}^{l_1'}
=0, \quad j= 1,\dots,m(l_0') 
$$
because
$\alpha_{N_0} \ne \beta_{N_0}$
or
$E_{i(x_{n_0})}^{l_1'}\cdot E_{i(y_{n_0})}^{l_1'} =0$.
This implies that
$$
E_{i(y_{N_0})}^{l_1'} 
S^*_{\beta_1 \cdots \beta_{N_0}} 
u_0 S_{\alpha_1 \cdots \alpha_{N_0}}
E_{i(x_{N_0})}^{l_1'}=0
$$
so that 
$$
E_{i(y_{N_0})}^{l_1'} 
S^*_{\beta_1 \cdots \beta_{N_0}} 
u S_{\alpha_1 \cdots \alpha_{N_0}}
E_{i(x_{N_0})}^{l_1'}=0
$$
 a contradiction to the preceding lemma. 
\qed
\enddemo
Thus we have
\proclaim{Lemma 4.5} 
For a partial isometry $u \in \FL$
satisfying 
$$
u \DL u^* \subset \DL, \qquad u^* \DL u \subset \DL,
$$
let
$\tau_u : \supp(u^* u) \rightarrow \supp(u u^* )$
be the homeomorphism
defined by
$Ad(u)(g) = g \circ {\tau_u^{-1}}$ for $ g \in \DL u^*u$.
Then
there exists $k_u \in {\Bbb N}$
such that
$$
\sigma_{\frak L}^{k_u}(\tau_u(x)) = \sigma_{\frak L}^{k_u}(x) \qquad \text{ for } 
x \in \supp(u^* u).
$$
\endproclaim
Therefore  by Lemma 4.2 and Lemma 4.5 we have
\proclaim{Proposition 4.6}
For any $v \in N_s(\OL,\DL)$,
the partial homomorphism $\tau_v$ 
induced by $Ad(v)$ on $\DL$
gives rise to an element of the topological full inverse semigroup
$[\sigma_{\frak L}]_{sc}$. 
If in particular $v$ belongs to
$N(\OL,\DL)$, 
then $\tau_v$ belongs to 
$[\sigma_{\frak L}]_{c}$.  
\endproclaim
\demo{Proof}
The argument of the proof is the same as
that of [Ma5;Proposition 4.7].
\qed
\enddemo
The unitaries ${\Cal U}(\DL)$ are naturally embedded into 
$N_s(\OL,\DL)$.
 We denote the embedding by $\id$.
For $v \in  N_s(\OL,\DL)$, 
the induced partial homemorphism $\tau_v$ 
on $X_{\frak L}$ 
gives rise to an element of $[\sigma_{\frak L}]_{sc}$
by the above proposition.
We then have
\proclaim{Theorem 4.7}
The diagrams  
$$
\CD
1 @ >>>  {\Cal U}(\DL) @>\id>> N(\OL,\DL)
  @>\tau>> [\sigma_{\frak L} ]_c @>>>1 \\
@. @|
@VV{\iota}V 
@VV{\iota}V @. \\
1 @ >>>  {\Cal U}(\DL) @>\id>> N_s(\OL,\DL)
  @>\tau>> [\sigma_{\frak L} ]_{sc} @>>>1. \\
\endCD
$$
are all commutative,
where two vertical arrows denoted by $\iota$ are inclusions.
The first row sequence is exact and splits as group,
and the second row sequence is exact and  splits as inverse semigroup.
\endproclaim
\demo{Proof}
By Proposition 4.6, the map
$\tau: v \in N_s(\OL,\DL)
\longrightarrow 
\tau_v \in [\sigma_{\frak L} ]_{sc}
$
defines a homomorphism as inverse semigroup such that
$\tau(N(\OL,\DL)) = [\sigma_{\frak L} ]_{c}$.
It is   surjective
by Proposition 4.1.
Suppose that $\tau_v =\id $ on $X_{\frak L}$ 
for some $v \in N_s(\OL,\DL)$.
This means that $Ad(v) =\id $
on
$\DL$.
Hence $v$ commutes 
with all of elements of $\DL$.
By Lemma 2.1, $v$ belongs to $\DL$.
Therefore the second row sequence is exact.
Similarly, the first row sequence is exact.
As in Proposition 4.1,
the partial isometry $u_\tau$ for $\tau \in [\sigma_{\frak L}]_{sc}$
defined by 
$u_\tau e_x = e_{\tau(x)}, \, x \in X_{\tau}\subset X_{\frak L}$ 
gives rise to sections of the both exact sequences.
Hence the both row sequences split.
The commutativity of the diagrams is clear
\qed
\enddemo
\heading 5. Orbit equivalence of $(X_{\frak L},\sigma_{\frak L})$
\endheading
In this section, 
we will study orbit equivalence between two dynamical systems
$(X_{{\frak L}_1}, \sigma_{{\frak L}_1})$ and
$(X_{{\frak L}_1}, \sigma_{{\frak L}_1})$
defined by $\lambda$-graph systems ${\frak L}_1$ and ${\frak L}_2$
respectively.

\noindent
{\bf Definition.}
For $\lambda$-graph systems
${\frak L}_1$ and ${\frak L}_2$,
if 
there exists a homeomorphism 
$h:X_{{\frak L}_1} \rightarrow X_{{\frak L}_2}$ 
such that 
$h(orb_{\sigma_{{\frak L}_1}}(x)) = orb_{\sigma_{{\frak L}_2}}(h(x))$ 
for $x \in X_{{\frak L}_1}$,
then  $(X_{{\frak L}_1}, \sigma_{{\frak L}_1})$ 
and $(X_{{\frak L}_2}, \sigma_{{\frak L}_2})$ 
are said to be topologically orbit equivalent.
In this case,  
there exist functions
$k_1,\, l_1:X_{{\frak L}_1} \rightarrow \Zp$
and
$k_2,\, l_2:X_{{\frak L}_2} \rightarrow \Zp$
satisfying
$$
\cases
\sigma_{{\frak L}_2}^{k_1(x)} (h(\sigma_{{\frak L}_1}(x))) 
=&  \sigma_{{\frak L}_2}^{l_1(x)}(h(x))\quad 
\text{ for } x \in X_{{\frak L}_1},\\
\sigma_{{\frak L}_1}^{k_2(y)} (h^{-1}(\sigma_{{\frak L}_2}(y))) 
=&  \sigma_{{\frak L}_1}^{l_2(y)}(h^{-1}(y))\quad 
\text{ for } y \in X_{{\frak L}_2}.
\endcases
\tag 5.1
$$

We say that 
$(X_{{\frak L}_1}, \sigma_{{\frak L}_1})$ 
and $(X_{{\frak L}_2}, \sigma_{{\frak L}_2})$ 
are  {\it continuously orbit equivalent}\,\
if there exist continuous functions 
$k_1,\, l_1:X_{{\frak L}_1} \rightarrow \Zp$
and
$k_2,\, l_2:X_{{\frak L}_2} \rightarrow \Zp$
satisfying the equalities 
(5.1).

The following lemma is straightforward.
\proclaim{Lemma 5.1}
If $h: X_{{\frak L}_1} \rightarrow X_{{\frak L}_2}$ is a homeomorphism satisfying
$\sigma_{{\frak L}_2}^{k(x)}(h(\sigma_{{\frak L}_1}(x))) 
= \sigma_{{\frak L}_2}^{l(x)}(h(x)), x \in X_{{\frak L}_1}$
for some functions
$k,l:X_{{\frak L}_1} \rightarrow \Zp$,
then by putting
$$
k^n(x) = \sum_{i=0}^{n-1}k(\sigma_{{\frak L}_1}^i(x)),\qquad
l^n(x) = \sum_{i=0}^{n-1}l(\sigma_{{\frak L}_1}^i(x)), \qquad n \in \Bbb N
$$
we have
$$
\sigma_{{\frak L}_2}^{k^n(x)}(h(\sigma_{{\frak L}_1}^n(x))) 
= \sigma_{{\frak L}_2}^{l^n(x)}(h(x)),
\qquad x \in X_{{\frak L}_1}.
$$
\endproclaim
\proclaim{Lemma 5.2}
If $h: X_{{\frak L}_1} \rightarrow X_{{\frak L}_2}$ is a homeomorphism satisfying (5.1),
then it satisfies
$$
h(orb_{\sigma_{{\frak L}_1}}(x)) = orb_{\sigma_{{\frak L}_2}}(h(x))\qquad 
\text{ for }x \in X_{{\frak L}_1}.
$$
Hence continuous orbit equivalence implies topological orbit equivalence.
\endproclaim
\demo{Proof}
By the preceding lemma,
one has
$$
h(\sigma_{{\frak L}_1}^n(x))
\subset
\sigma_{{\frak L}_2}^{-k^n(x)}( 
 \sigma_{{\frak L}_2}^{l^n(x)}(h(x))),
\qquad x  \in X_{{\frak L}_1}, n \in \Bbb N
$$
so that 
$h(\sigma_{{\frak L}_1}^n(x))
\subset
orb_{\sigma_{{\frak L}_2}}(h(x)).
$
For 
$(z_1,\dots,z_m,x_1,x_2,\dots ) \in \sigma_{{\frak L}_1}^{-m}(x)$,
where
$x = (x_n)_{n \in \Bbb N},$
one has
$
\sigma^m(z_1,\dots,z_m,x_1,x_2,\dots )= x
$
and hence
$h(z_1,\dots,z_m,x_1,x_2,\dots ) \in 
\sigma_{{\frak L}_2}^{-l_1^m(x)}
\sigma_{{\frak L}_2}^{-k_1^m(x)}(h(x)).
$ 
This implies that
$
h(orb_{\sigma_{{\frak L}_1}}(x)) \subset orb_{\sigma_{{\frak L}_2}}(h(x)).
$

One similarly has the inclusion relation 
$
h^{-1}(orb_{\sigma_{{\frak L}_2}}(y)) \subset orb_{\sigma_{{\frak L}_1}}(h^{-1}(y))
$ for $ y \in X_{{\frak L}_2}$
by considering $h^{-1}$ as $h $ 
in the above discussion.
This implies that
$
orb_{\sigma_{{\frak L}_2}}(h(x)) \subset h( orb_{\sigma_{{\frak L}_1}}(x))
$ for $ x \in X_{{\frak L}_1}$
so that
$
h(orb_{\sigma_{{\frak L}_1}}(x)) = orb_{\sigma_{{\frak L}_2}}(h(x)).
$ 
\qed
\enddemo

\proclaim{Proposition 5.3}
If there exists a homeomorphism
$h:X_{{\frak L}_1} \longrightarrow X_{{\frak L}_2}$
such that
$ h \circ [\sigma_{{\frak L}_1}]_{sc} \circ h^{-1} = [\sigma_{{\frak L}_2}]_{sc}$,
then 
$(X_{{\frak L}_1}, \sigma_{{\frak L}_1})$ 
and  
$(X_{{\frak L}_2}, \sigma_{{\frak L}_2})$
are continuously orbit equivalent.
\endproclaim
\demo{Proof}
Let us denote by
$\{v_1^2,\dots,v_{m(2)}^2\}$
 the vertex set
 $V_2$. 
 For $i=1,\dots,m(2)$, 
 let 
$B_2(v_i^2)$ 
be the set of all admissible words of length $2$
terminating at $v_i^2$.
That is
$$
\align
B_2(v_i^2) = \{ (\mu_1,\mu_2) \in B_2(X_\Lambda) \mid 
& \text{there exist } e_1\in E_{0,1}, e_2\in E_{1,2}; \\
 \lambda(e_1) = \mu_1, & \lambda(e_2) = \mu_2, 
 t(e_1) = s(e_2), t(e_2) = v_i^2 \}. 
\endalign
$$
For $\mu \in B_2(v_i^2)$,
by Lemma 3.1, 
there exists 
$\tau_\mu \in [\sigma_{{\frak L}_1}]_{sc}$ 
such that
$
\tau_\mu(x) = \sigma_{\frak L}(x)$ for $x \in U_{\mu,v_i^2}$.
Put
$
\tau_{h,\mu} = h \circ \tau_\mu \circ h^{-1} 
\in 
h\circ [\sigma_{{\frak L}_1}]_{sc} \circ h^{-1} = [\sigma_{{\frak L}_2}]_{sc}.
$
There exist
$k_{\tau_{h,\mu}}, l_{\tau_{h,\mu}}: 
h(U_{\mu,v_i^2}) \rightarrow \Zp
$
such that
$$
\sigma_{{\frak L}_2}^{k_{\tau_{h,\mu}}(y)}(\tau_{h,\mu}(y)) = 
\sigma_{{\frak L}_2}^{l_{\tau_{h,\mu}}(y)}(y), \qquad
y \in h(U_{\mu,v_i^2}).
$$
For $x\in U_{\mu,v_i^2}$, 
one has 
$
\tau_{h,\mu}(h(x)) 
= h \circ \tau_\mu(x) = h \circ \sigma_{{\frak L}_1}(x)
$
so that 
$$
\sigma_{{\frak L}_2}^{k_{\tau_{h,\mu}}(h(x))}(h \circ \sigma_{{\frak L}_1}(x)) 
= 
\sigma_{{\frak L}_2}^{l_{\tau_{h,\mu}}(h(x))}(h(x)),
\qquad
x \in U_{\mu,v_i^2}.
$$
Since 
$X_{{\frak L}_1} $
is a disjoint union
$
\cup_{i=1}^{m(2)} \cup_{\mu \in B_2(v_i^2)} U_{\mu,v_i^2},
$
by putting
$$
k_1(x) = k_{\tau_{h,\mu}}(h(x)),
\quad
l_1(x) = l_{\tau_{h,\mu}}(h(x))
\qquad
\text{ for }
x \in  U_{\mu,v_i^2},
$$
we have
continuous functions
$k_1, l_1 : X_{{\frak L}_1}\longrightarrow \Zp$
satisfying
$$
\sigma_{{\frak L}_2}^{k_1(x)}(h \circ \sigma_{{\frak L}_1}(x)) 
= 
\sigma_{{\frak L}_2}^{l_1(x)}(h(x)),
\quad
x \in X_{{\frak L}_1}.
$$
We similarly have 
continuous functions
$k_2, l_2 : X_{{\frak L}_2}\longrightarrow \Zp$
satisfying
$$
\sigma_{{\frak L}_1}^{k_2(y)}(h^{-1} \circ \sigma_{{\frak L}_2}(y)) 
= 
\sigma_{{\frak L}_1}^{l_2(x)}(h^{-1}(y)),
\quad
y \in X_{{\frak L}_2}.
$$
Hence
$(X_{{\frak L}_1}, \sigma_{{\frak L}_1})$ 
and  
$(X_{{\frak L}_2}, \sigma_{{\frak L}_2})$
are continuously orbit equivalent.
\qed
\enddemo 
Conversely we have
\proclaim{Proposition 5.4}
If 
$(X_{{\frak L}_1}, \sigma_{{\frak L}_1})$ 
and  
$(X_{{\frak L}_2}, \sigma_{{\frak L}_2})$
are continuously orbit equivalent,
then there exists a homeomorphism
$h:X_{{\frak L}_1} \longrightarrow X_{{\frak L}_2}$
such that
$ h \circ [\sigma_{{\frak L}_1}]_{sc} \circ h^{-1} 
= [\sigma_{{\frak L}_2}]_{sc}$.
\endproclaim
\demo{Proof}
Suppose that 
there exist a homeomorphism 
$h:X_{{\frak L}_1} \rightarrow X_{{\frak L}_2}$ 
and
continuous functions
$k_1,\, l_1:X_{{\frak L}_1} \rightarrow \Zp$
and
$k_2,\, l_2:X_{{\frak L}_2} \rightarrow \Zp$
satisfying (5.1).
For $n \in {\Bbb N}$,
let 
$k_1^n, l_1^n:X_{{\frak L}_1} \longrightarrow \Zp$
and 
$k_2^n, l_2^n:X_{{\frak L}_2} \longrightarrow \Zp$
be continuous functions as in Lemma 5.1 such that
$$
\sigma_{{\frak L}_2}^{k_1^n(x)} (h(\sigma_{{\frak L}_1}^n(x)) 
= \sigma_{{\frak L}_2}^{l_1^n(x)}(h(x)),
\qquad
\sigma_{{\frak L}_1}^{k_2^n(y)} (h^{-1}(\sigma_{{\frak L}_2}^n(y)) 
= \sigma_{{\frak L}_1}^{l_2^n(y)}(h^{-1}(y)) 
\tag 5.2
$$
for $x \in X_{{\frak L}_1}$ and $y \in X_{{\frak L}_2}$.
For any $\tau \in [\sigma_{{\frak L}_1}]_{sc}$, 
there exist continuous functions:
$k_\tau, l_\tau : X_{{\frak L}_1} \longrightarrow \Zp$
such that
$$
\sigma_{{\frak L}_1}^{k_\tau(x)}(\tau(x)) 
= \sigma_{{\frak L}_1}^{l_\tau(x)}(x),
\qquad
x \in X_{\tau}. \tag 5.3
$$
For $ y \in h(X_\tau)$,
set
$x = h^{-1}(y) \in X_\tau$.
Put
$m=k_\tau(x)$.
By (5.2) and (5.3), one has
$$
\sigma_{{\frak L}_2}^{l_1^m (\tau(x))}(h(\tau(x))
=\sigma_{{\frak L}_2}^{k_1^m(\tau(x))} (h(\sigma_{{\frak L}_1}^m(\tau(x)))
=\sigma_{{\frak L}_2}^{k_1^m(\tau(x))} (h(\sigma_{{\frak L}_1}^{l_\tau(x)}(x))
$$
Put
$n=l_\tau(x)\in \Bbb N$.
By applying $\sigma_{{\frak L}_2}^{k_1^n(x)}$ 
to the above equality,
one has by (5.2)
$$
\align
& \sigma_{{\frak L}_2}^{k_1^n(x) +l_1^m (\tau(x))}(h(\tau(x))\\
=& \sigma_{{\frak L}_2}^{k_1^m(\tau(x))}\sigma_{{\frak L}_2}^{k_1^n(x)} 
(h(\sigma_{{\frak L}_1}^n(x)))
=  \sigma_{{\frak L}_2}^{k_1^m(\tau(x))}\sigma_{{\frak L}_2}^{l_1^n(x)} (h(x))
=  \sigma_{{\frak L}_2}^{k_1^m(\tau(x))+ l_1^n(x)} (h(x))
\endalign
$$
and hence
$$
\sigma_{{\frak L}_2}^{k_1^n(x) +l_1^m (\tau(x))}(h\circ \tau \circ h^{-1}(y))
=\sigma_{{\frak L}_2}^{k_1^m(\tau(x))+ l_1^n(x)} (y).
$$
By setting for $y \in h(X_\tau)$,
$$
\align
k_\tau^h(y) 
& = k_1^n(x) +l_1^m (\tau(x))
= k_1^{l_\tau( h^{-1}(y))}(h^{-1}(y)) 
+l_1^{k_\tau(h^{-1}(y))} (\tau(h^{-1}(y))),\\
l_\tau^h(y) 
& = k_1^m(\tau(x)) +l_1^n (x)
= k_1^{k_\tau( h^{-1}(y))}(\tau(h^{-1}(y))) 
+l_1^{l_\tau(h^{-1}(y))} (h^{-1}(y)),
\endalign
$$
one has
$$
\sigma_{{\frak L}_2}^{k_\tau^h(y)}(h\circ \tau \circ h^{-1}(y))
=\sigma_{{\frak L}_2}^{l_\tau^h(y)} (y)
\qquad 
\text{ for  }
y \in h(X_{\tau})
$$
so that
$
h\circ \tau \circ h^{-1} \in [\sigma_{{\frak L}_2}]_{sc}
$
and hence
$
h \circ [\sigma_{{\frak L}_1}]_{sc} \circ h^{-1} 
\subset [\sigma_{{\frak L}_2}]_{sc}.
$
Similarly one has
$
h^{-1}\circ  [\sigma_{{\frak L}_2}]_{sc} \circ h 
\subset [\sigma_{{\frak L}_1}]_{sc}
$
and concludes 
$
h \circ [\sigma_{{\frak L}_1}]_{sc} \circ h^{-1} = [\sigma_{{\frak L}_2}]_{sc}.
$
\qed
\enddemo

\proclaim{Proposition 5.5}
If there exists an isomorphism 
$\Psi: {\Cal O}_{{\frak L}_1} \longrightarrow {\Cal O}_{{\frak L}_2}$
such that
$\Psi(\DLA) = \DLB$,
then there exists a homeomorphism
$h:X_{{\frak L}_1} \longrightarrow X_{{\frak L}_2}$
such that
$ h \circ [\sigma_{{\frak L}_1}]_{sc} \circ h^{-1} 
= [\sigma_{{\frak L}_2}]_{sc}$.
\endproclaim
\demo{Proof}
Suppose that 
there exists an isomorphism 
$\Psi: {\Cal O}_{{\frak L}_1} \longrightarrow {\Cal O}_{{\frak L}_2}$
such that
$\Psi(\DLA) = \DLB$.
By the split exact sequences 
$$
1
\longrightarrow
{\Cal U}({\Cal D}_{{\frak L}_i})
\longrightarrow
N_s({\Cal O}_{{\frak L}_i},{\Cal D}_{{\frak L}_i})
\longrightarrow
[\sigma_{{\frak L}_i}]_{sc}
\longrightarrow
1,
\qquad i=1,2
$$
of
inverse semigroups,
one may find an isomorphism
$\widetilde{\Psi}: 
[\sigma_{{\frak L}_1}]_{sc} \longrightarrow [\sigma_{{\frak L}_2}]_{sc}$
of inverse semigroup
such that the following diagrams are commutative:
$$
\CD
1 @ >>>  {\Cal U}(\DLA) @>\id>> N_s(\OLA,\DLA)
  @>\tau>> [\sigma_{{\frak L}_1} ]_{sc} @>>>1 \\
@. @VV{\Psi |_{{\Cal U}(\DLA)}}V
@VV{\Psi}V 
@VV{\widetilde{\Psi}}V @. \\
1 @ >>>  {\Cal U}(\DLB) @>\id>> N_s(\OLB,\DLB)
  @>\tau>> [\sigma_{{\frak L}_2} ]_{sc} @>>>1. \\
\endCD
$$
For $v \in N_s(\OLA,\DLA)$,
take the partial homeomorphism 
$\tau_v: X_v \longrightarrow Y_v$ 
satisfying
$Ad(v)(f) = f \circ \tau_v^{-1}$ for $f \in C(X_v)$. 
Let $h:X_{{\frak L}_1} \longrightarrow X_{{\frak L}_2}$ 
be the homeomorphism
satisfying
$\Psi(f) = f \circ h^{-1}$ for $f \in C(X_{{\frak L}_1})$.
For $g \in C(h(X_v))$,
we have
$$
\Psi \circ Ad(v) \circ \Psi^{-1} (g)
=
g \circ h \circ \tau_v^{-1}\circ h,
\quad
\text{ and }
\quad
Ad(\Psi(v))(g) =g \circ \tau_{\Psi(v)}^{-1}.
$$
By
the identity 
$
\Psi \circ Ad(v) \circ \Psi^{-1} = Ad(\Psi(v)),
$
one has
$$
g \circ h \circ \tau_v^{-1}\circ h
=
g \circ \tau_{\Psi(v))}^{-1}
\quad
\text{ for }
g \in C(h(X_v)).
$$
Hence
$h \circ \tau_v \circ h^{-1}
=
\tau_{\Psi(v)}.
$
As
$
[\sigma_{{\frak L}_i}]_{sc}
=\{ \tau_v \mid v \in 
N_s({\Cal O}_{{\frak L}_i},{\Cal D}_{{\frak L}_i})\}, i=1,2,
$
one sees that
$h \circ [\sigma_{{\frak L}_1}]_{sc} \circ h^{-1} 
= [\sigma_{{\frak L}_2}]_{sc}.$
\qed
\enddemo

\proclaim{Proposition 5.6}
If $(X_{{\frak L}_1}, \sigma_{{\frak L}_1})$ 
and  
$(X_{{\frak L}_2}, \sigma_{{\frak L}_2})$
are continuously orbit equivalent,
then 
there exists an isomorphism $\Psi:\OLA \longrightarrow \OLB$
such that
$\Psi(\DLA) = \DLB$.
\endproclaim
\demo{Proof}
The proof is essentially same as the proof of 
Proposition 4.1 and [Ma5:Proposition 5.5].
We omit its proof.
\qed
\enddemo
Therefore we have
\proclaim{Theorem 5.7}
Let ${\frak L}_1$ and ${\frak L}_2$ be 
$\lambda$-graph systems satisfying condition (I).
The following  are  equivalent:
 \roster
 \item There exists an isomorphism 
 $\Psi: \OLA \rightarrow  \OLB$
 such that $\Psi(\DLA) = \DLB$.
\item
$(X_{{\frak L}_1}, \sigma_{{\frak L}_1})$ 
and $(X_{{\frak L}_2}, \sigma_{{\frak L}_2})$ 
are continuously orbit equivalent.
\item 
There exists a homeomorphism 
$h: X_{{\frak L}_1} \rightarrow X_{{\frak L}_2}$ such that
$h \circ [\sigma_{{\frak L}_1} ]_{sc} \circ h^{-1} 
= [\sigma_{{\frak L}_2} ]_{sc}$. 
\endroster
\endproclaim

\noindent
{\bf Example.}

Let $G=(V,E)$ 
be a finite directed graph with 
$V = \{ v_1,v_2 \}$ and 
$E = \{ e, f, g \}$
such that
$$
s(e) = t(e) = s(f) = t(g) = v_1,\qquad 
t(f) = s(g) = v_2.
$$
Put the alphabet sets 
 $\Sigma_1 = \{ \bold{1}, \bold{2} \}$
and
$\Sigma_2 = \{ \alpha, \beta \}$.
Define two labeling maps
$\lambda_i: E \longrightarrow \Sigma_i, i=1,2$
by setting
$$
\lambda_1(e) = \lambda_1(f) = \bold{1}, \enskip \lambda_1(g) = \bold{2},\qquad
\lambda_2(e) = \alpha, \enskip \lambda_2(f) = \lambda_2(g) = \beta.
$$
Let us denote by
${\Cal G}_i$  
the labeled graph
$(G, \lambda_i)$ over $\Sigma_i$ 
for $i=1,2$.
Hence their underlying directed graphs are both $G$.
The labeled graphs
${\Cal G}_1$  and 
${\Cal G}_2$  
have its adjacency matrices as
$$
\bmatrix
\bold{1} & \bold{1} \\
\bold{2} & 0
\endbmatrix,
\qquad
\bmatrix
\alpha & \beta \\
\beta & 0
\endbmatrix
$$
respectively.
Let ${\frak L}_i =(V^{(i)}, E^{(i)}, \lambda^{(i)}, \Sigma_i)$ 
be the $\lambda$-graph systems 
associated to the labeled graphs ${\Cal G}_i$ for $i=1,2$
respectively.
They are defined by setting
$$
V^{(i)}_{l,l+1} = V, \quad
E^{(i)}_{l,l+1} = E, \quad 
\lambda^{(i)} = \lambda_i
$$
for all $l \in \Zp$ and $i=1,2$.
We then have
$\Omega_{{\frak L}_i} = V = \{v_1,v_2\}, i=1,2$.
The correspondence:
$$
(\bold{1},v_1) \rightarrow (\alpha,v_1),\quad
(\bold{1},v_2) \rightarrow (\beta,v_2),\quad
(\bold{2},v_1) \rightarrow (\beta,v_1)       
$$
yields a homeomorphism
$h: X_{{\frak L}_1} \longrightarrow  X_{{\frak L}_2}$
that gives rise to a continuous orbit equivalence
between
$(X_{{\frak L}_1},\sigma_{{\frak L}_1})$
and
$(X_{{\frak L}_2},\sigma_{{\frak L}_2})$.
One indeed sees that the $C^*$-algebras $\OLA$ and $\OLB$
are both isomorphic to the Cuntz-Krieger algebra
${\Cal O}_F$
where
$F
=\left[\smallmatrix
1 & 1 \\
1 & 0
\endmatrix\right]
$,
although the subshift presented by the $\lambda$-graph system
${\frak L}_2$ is the even shift that is not a Markov shift.

\heading 6. Orbit equivalence of the factor map
$\pi^{\frak L}_\Lambda: X_{\frak L} \longrightarrow X_{\Lambda}$
\endheading
For a $\lambda$-graph system $\frak L$ over $\Sigma$,
let 
$\Lambda$ be the subshift presented by ${\frak L}$.
Then we have a factor map
$
\pi^{\frak L}_\Lambda: (X_{\frak L}, \sigma_{\frak L}) \longrightarrow 
(X_\Lambda,\sigma_\Lambda).
$
In this section, we will study orbit structure between 
two dynamical systems
$(X_{\frak L}, \sigma_{\frak L}) 
$
and
$ 
(X_\Lambda,\sigma_\Lambda)
$
through the factor map
$\pi^{\frak L}_\Lambda$.
\proclaim{Lemma 6.1}
$\pi^{\frak L}_\Lambda(orb_{\sigma_{\frak L}}(x)) 
= orb_{\sigma_\Lambda}(\pi^{\frak L}_\Lambda(x))$ for 
$ x\in X_{\frak L}$.
\endproclaim
\demo{Proof}
Take an arbitrary element $x = (x_n)_{n \in \Bbb N} \in X_{\frak L}$.
For $w \in orb_{\sigma_{\frak L}}(x)$,
we have
$w = (z_1,\dots,z_k,x_{l+1},x_{l+2},\dots ) \in X_{\frak L}$
for some
$z_1\cdots z_k \in B_k(X_{\frak L})$
and $l \in \Zp$.
It is easy to see that 
$$
\pi^{\frak L}_\Lambda(w)
\in 
\sigma_\Lambda^{-k}(\sigma_\Lambda^l(\pi^{\frak L}_\Lambda(x))) 
\subset 
orb_{\sigma_\Lambda}(\pi^{\frak L}_\Lambda(x)).
$$ 
Conversely,
put
$(\alpha_n)_{n \in \Bbb N} = \pi^{\frak L}_\Lambda(x)$.
Each element 
$a \in orb_{\sigma_\Lambda}(\pi^{\frak L}_\Lambda(x))$
has
of the form
$ a 
= (\gamma_1, \dots, \gamma_k,\alpha_{l+1},\alpha_{l+2},\dots ) 
\in X_\Lambda
$
for some $\gamma_1\cdots \gamma_k \in B_k(X_\Lambda)$
and
$l \in \Zp$.
Put
$v_0 = v_0(\sigma_{\frak L}^l(x)) \in \Omega_{\frak L}$.
Since $\frak L$ is left-resolving,
there uniquely exists
$v_{-1} \in \Omega_{\frak L}$
such that
$(v_{-1},\gamma_k,v_0) \in E_{\frak L}$.
Inductively 
there uniquely exist
$v_{-2}, v_{-3},\dots, v_{-k} \in \Omega_{\frak L}$
such that
$(v_{-i}, \gamma_{k-(i-1)}, v_{-(i-1)}) \in E_{\frak L}$
for $i=1,2,\dots, k$.
Put $z_{k-(i-1)} = (\gamma_{k-(i-1)}, v_{-(i-1)})$
for
$i=1,2,\dots, k$
so that
$w = (z_1,\dots,z_k,x_{l+1},x_{l+2},\dots ) \in X_{\frak L}$
and
$\pi^{\frak L}_\Lambda(w) = a$.
Since
$w \in \sigma_{\frak L}^{-k}(\sigma_{\frak L}^l(x)) 
\subset orb_{\sigma_{\frak L}}(x),
$
one has
$ a \in 
  \pi^{\frak L}_\Lambda(orb_{\sigma_{\frak L}}(x)).
$  
\qed
\enddemo
For $\lambda$-graph systems ${\frak L}_1$ and ${\frak L}_2$,
let
$\Lambda_1$ and $\Lambda_2$ be the subshifts presented by 
${\frak L}_1$ and ${\frak L}_2$ respectively.

\noindent
{\bf Definition.}
Two factor maps
$\pi^{{\frak L}_1}_{\Lambda_1}$
and
$\pi^{{\frak L}_2}_{\Lambda_2}$
are said to be continuously orbit equivalent
if there exist homeomorphisms
$h_{\frak L}:X_{{\frak L}_1}\longrightarrow X_{{\frak L}_1}$
and
$h_{\Lambda}:X_{\Lambda_1}\longrightarrow X_{\Lambda_2}$
such that
$
\pi^{{\frak L}_2}_{\Lambda_2} \circ h_{\frak L} 
= h_\Lambda \circ \pi^{{\frak L}_1}_{\Lambda_1}
$
and
continuous functions
$k_1,l_1:X_{{\frak L}_1}\longrightarrow \Zp$
and
$k_2,l_2:X_{{\frak L}_2}\longrightarrow \Zp$
such that 
$$
\align
\sigma_{{\frak L}_2}^{k_1(x)}(h_{\frak L} \circ \sigma_{{\frak L}_1}(x)) 
& = 
\sigma_{{\frak L}_2}^{l_1(x)}(h_{\frak L}(x)),
\quad
x \in X_{{\frak L}_1} \tag 6.1\\
\sigma_{{\frak L}_1}^{k_2(y)}(h_{\frak L}^{-1} \circ \sigma_{{\frak L}_2}(y)) 
& = 
\sigma_{{\frak L}_1}^{l_2(x)}(h_{\frak L}^{-1}(y)),
\quad
y \in X_{{\frak L}_2}.\tag 6.2
\endalign
$$
We note that 
the equalities (6.1) and (6.2) imply 
$$
h_{\frak L}(orb_{\sigma_{{\frak L}_1}}(x)) = 
orb_{\sigma_{{\frak L}_2}}(h_{\frak L}(x))
\qquad
\text{ for } x \in X_{{\frak L}_1}.\tag 6.3
$$
\proclaim{Lemma 6.2}
Suppose that two factor maps
$\pi^{{\frak L}_1}_{\Lambda_1}$ and $\pi^{{\frak L}_2}_{\Lambda_2}$
are continuously orbit equivalent and keep the above notation. 
Then we have
\roster
\item"(i)"
$$
\align
\sigma_{\Lambda_2}^{k_1(x)}(h_\Lambda \circ \sigma_{\Lambda_1}
(\pi^{{\frak L}_1}_{\Lambda_1}(x)) 
& = 
\sigma_{\Lambda_2}^{l_1(x)}(h_\Lambda(\pi^{{\frak L}_1}_{\Lambda_1}(x)),
\quad
x \in X_{{\frak L}_1}, \\
\sigma_{\Lambda_1}^{k_2(y)}(h_\Lambda^{-1} \circ \sigma_{\Lambda_2}
(\pi^{{\frak L}_2}_{\Lambda_2}(y)) 
& = 
\sigma_{\Lambda_1}^{l_2(y)}(h_\Lambda^{-1}(\pi^{{\frak L}_2}_{\Lambda_2}(y)),
\quad
y \in X_{{\frak L}_2}.
\endalign
$$
\item"(ii)"
$$
h_\Lambda(orb_{\sigma_{\Lambda_1}}(a)) = 
orb_{\sigma_{\Lambda_2}}(h_\Lambda(a))
\qquad
\text{ for } a \in X_{\Lambda_1}.
$$
\endroster
\endproclaim
\demo{Proof}
(i) follows from (6.1) and (6.2),
and
(ii) follows from (6.3).
\qed
\enddemo
The following lemma is direct.
\proclaim{Lemma 6.3}
Two factor maps
$\pi^{{\frak L}_1}_{\Lambda_1}$
and
$\pi^{{\frak L}_2}_{\Lambda_2}$
are  continuously orbit equivalent
if and only if
there exists a homeomorphism
$h_{\frak L}:X_{{\frak L}_1}\longrightarrow X_{{\frak L}_2}$
that yields a continuously orbit equivalence between
$(X_{{\frak L}_1}, \sigma_{{\frak L}_1})$ 
and $(X_{{\frak L}_2}, \sigma_{{\frak L}_2})$
 and there exists a homemorphism
$h_{\Lambda}:X_{\Lambda_1}\longrightarrow X_{\Lambda_2}$
such that
$
\pi^{{\frak L}_2}_{\Lambda_2} \circ h_{\frak L} 
= h_\Lambda \circ \pi^{{\frak L}_1}_{\Lambda_1}.
$
\endproclaim
We note that the factor map
$\pi^{\frak L}_\Lambda:X_{\frak L} \longrightarrow X_\Lambda$
induces an embedding of
$C(X_\Lambda)$ into
$C(X_{\frak L})$,
that corresponds to the natural embedding
of $\DLambda$ into $\DL$. 
Let
$N_s(\OL,\DLambda)$be the set of all partial isometries 
$v \in \OL$ such that
$v \DLambda v^* \subset \DLambda$
and
$v^* \DLambda v \subset \DLambda$.
\proclaim{Lemma 6.4}
$N_s(\OL,\DLambda) \subset N_s(\OL,\DL).$
\endproclaim
\demo{Proof}
For $v \in N_s(\OL,\DLambda)$, and
$x \in \DL, a \in \DLambda$, we have
$$
v x v^* a = v x v^* a v v^* = v v^* a v x v^* = a v x v^*
$$
so that 
$v x v^* \in \DLambda' \cap \OL = \DL$.
Hence 
$v \DL v^* \subset \DL$,
 and similarly
 $v^* \DL v \subset \DL$.
 This implies that 
 $v \in N_s(\OL,\DL).$
\qed
\enddemo
Suppose that both $\lambda$-graph systems
${\frak L}_1$ and 
${\frak L}_2$
satisfy condition (I).
\proclaim{Lemma 6.5}
If there exists an isomorphism
$\Psi:\OLA \longrightarrow \OLB$ such that
$\Psi(\DLambdaA)=\DLambdaB$,
then
$\Psi(\DLA) =\DLB$. 
\endproclaim
\demo{Proof}
Suppose that
$\Psi(\DLambdaA)=\DLambdaB$.
For $x\in \DLA$ and $b \in \DLambdaB$,
take $a \in \DLambdaA$ such that $\Psi(a) =b$.
It then follows that
$$
\Psi(x) b = \Psi(xa) = \Psi(a)\Psi(x) = b \Psi(x)
$$
so that
$\Psi(x) $ commutes with all elements of $\DLambdaB$,
and hence
$\Psi(x) \in \DLB$.
This implies that  
$\Psi(\DLA) \subset \DLB$.
Similarly we have
 $\Psi^{-1}(\DLB) \subset \DLA$
 so that $\Psi(\DLA)= \DLB$.
\qed
\enddemo
\proclaim{Theorem 6.6}
Let ${\frak L}_1$ and
${\frak L}_2$ be $\lambda$-graph systems
satisfying condition (I).
Let
$X_{\Lambda_1}$ and $X_{\Lambda_2}$ be their respect right one-sided subshifts.
The following are equivalent:
\roster
\item"(i)"
There exists an isomorphism
$\Psi:\OLA \longrightarrow \OLB$ such that
$\Psi(\DLambdaA)=\DLambdaB$.
\item"(ii)"
The factor maps $\pi^{{\frak L}_1}_{\Lambda_1}$
and
$\pi^{{\frak L}_2}_{\Lambda_2}$
are
continuously orbit equivalent.
\item"(iii)" 
There exist homeomorphisms
$h_{\frak L}: X_{{\frak L}_1} \longrightarrow
 X_{{\frak L}_2}$
and
$h_\Lambda:X_{\Lambda_1} \longrightarrow X_{\Lambda_2}$
 such that
$
\pi^{{\frak L}_2}_{\Lambda_2} \circ h_{\frak L} 
= h_\Lambda \circ \pi^{{\frak L}_1}_{\Lambda_1}$
and
 $h_{\frak L}\circ [\sigma_{{\frak L}_1}]_{sc}\circ h_{\frak L}^{-1}
 = [\sigma_{{\frak L}_2}]_{sc}.
$
\endroster
\endproclaim
\demo{Proof}
(ii)$ \Leftrightarrow$ (iii):
The equivalence between (ii) and (iii)  comes from Lemma 6.3.

(i)$ \Rightarrow $(iii):
Suppose that
there exists an isomorphism
$\Psi:\OLA \longrightarrow \OLB$ such that
$\Psi(\DLambdaA)=\DLambdaB$.
By Lemma 6.5,
one has
$\Psi(\DLA)=\DLB$.
Let $h_{\frak L}:X_{{\frak L}_1}\rightarrow X_{{\frak L}_2}$
be the homeomorphism
induced by 
$\Psi:\DLA\longrightarrow \DLB$
such that
$\Psi(f) = f \circ h^{-1}$
for
$f \in \DLA$.
Then $h_{\frak L}$ satisfies
$
h \circ [\sigma_{{\frak L}_1}]_{sc} \circ h^{-1} 
= [\sigma_{{\frak L}_2}]_{sc}
$
by Proposition 5.5.
Since
$\Psi(\DLambdaA)=\DLambdaB$,
there exists a homeomorphism
$h_{\Lambda}: X_{\Lambda_1} \longrightarrow X_{\Lambda_2}
$
such that
$
h_\Lambda\circ\pi^{{\frak L}_1}_{\Lambda_1} 
= \pi^{{\frak L}_2}_{\Lambda_2}\circ h_{\frak L}$.  

(ii)$ \Rightarrow $(i):
Suppose that the
factor maps $\pi^{{\frak L}_1}_{\Lambda_1}$
and
$\pi^{{\frak L}_2}_{\Lambda_2}$
are
continuously orbit equivalent.
Since
$(X_{{\frak L}_1},\sigma_{{\frak L}_1})$
and
$(X_{{\frak L}_2},\sigma_{{\frak L}_2})$
are continuously orbit equivalent,
by Proposition 5.6
there exists an isomorphism
$\Psi:\OLA \longrightarrow \OLB$
such that
$\Psi(\DLA)=\DLB$
and
$\Psi(f) = f \circ h_{\frak L}^{-1}$ 
for
$f \in \DLA$.
For $g \in \DLambdaA$,
one sees that
$g \circ \pi^{{\frak L}_1}_{\Lambda_1} \in \DLA$
so that
$$
\Psi(g \circ \pi^{{\frak L}_1}_{\Lambda_1})
=g \circ \pi^{{\frak L}_1}_{\Lambda_1}\circ  h_{\frak L}^{-1}
= g \circ h_\Lambda^{-1}\circ \pi^{{\frak L}_2}_{\Lambda_2}
$$
This means that
$\Psi(\DLambdaA) \subset\DLambdaB$,
and similarly
$\Psi^{-1}(\DLambdaB) \subset\DLambdaA$.
Therefore we conclude that
$\Psi(\DLambdaA)=\DLambdaB.$
\qed
\enddemo


\heading 7. Orbit equivalence of one-sided subshifts
\endheading

Let $\Lambda$ be a two-sided subshift over $\Sigma$
and $X_\Lambda$ its right one-sided subshift.
The canonical $\lambda$-graph system 
${\frak L}^\Lambda$ for $\Lambda$
is defined as in the following way ([Ma2]).
For $a =(a_n)_{n \in \Bbb N} \in X_\Lambda$
and $l \in \Zp$,
denote by $P_l(a)$ the predecessor set of length $l$ of $a$, 
that is
$$
P_l(a) = \{ (\mu_1,\dots,\mu_l) \in B_l(X_\Lambda)\mid
(\mu_1,\dots,\mu_l,a_1,a_2,\dots ) \in X_\Lambda \}.
$$
Two sequences 
$a = (a_n )_{n \in \Bbb N}$ 
and
$b = (b_n )_{n \in \Bbb N}$ 
in $X_{\Lambda}$
are said to be $l$-past equivalent if $P_l(a) = P_l(b)$,
and written as
$a\underset{l}\to{\sim}b$.
The equivalence class of $a$ in $X_\Lambda / \underset{l}\to{\sim}$
is denoted by $[a]_l$.
The vertex set $V_l$ of the $\lambda$-graph system 
is
the set
$X_\Lambda / \underset{l}\to{\sim}$.
We set
$v^l(a) = [a]_l$.
Then
$ (v^l(a))_{l \in \Zp}$
defines an $\iota$-orbit of
$\Omega_{{\frak L}^\Lambda}$,
denoted by
$v(a). $
An edge labeled $\alpha$ from 
$v^l(a)$ to $v^{l+1}(b)$ is defined if
$ a \underset{l}\to{\sim} (\alpha,b_1,b_2,\dots )$,
where
$b = (b_n)_{n \in \Bbb N}$.
\proclaim{Lemma 7.1}
For $a =(a_n)_{n \in \Bbb N} \in X_\Lambda$,
$(a_n, v_n(a))_{n\in \Bbb N}$ defines an element of $X_{{\frak L}^\Lambda}$.
\endproclaim
\demo{Proof}
For each $n \in \Bbb N$ and $l \in \Zp$, 
there is a unique  edge from 
$[(a_n,a_{n+1},\dots ) ]_l\in V_l $ to 
$[(a_{n+1},a_{n+2}, \dots ) ]_{l+1}\in V_{l+1} $
labeled $a_n$.
Hence
$(v_{n-1}(a), a_n, v_n(a))$ 
belongs to 
$E_{{\frak L}^{\Lambda}}$
for all $n \in \Bbb N$,
so that 
$(a_n, v_n(a))_{n\in \Bbb N}$ defines an element of $X_{{\frak L}^\Lambda}$.
\qed
\enddemo
We put the embedding of $X_\Lambda$ into $X_{{\frak L}^\Lambda}$:
$$
\iota_\Lambda: 
a =(a_n)_{n \in \Bbb N} \in X_\Lambda \longrightarrow
(a_n, v_n(a))_{n\in \Bbb N} \in X_{{\frak L}^\Lambda}.
$$
It is straightforward to see that the following lemma holds:
\proclaim{Lemma 7.2}
The map 
$
\iota_\Lambda: 
X_\Lambda \longrightarrow
X_{{\frak L}^\Lambda}
$
is injective and 
$
\iota_\Lambda(X_\Lambda)
$ is dense in $X_{{\frak L}^\Lambda}$.
\endproclaim
We endow $X_\Lambda$ with a new topology induced by the injection
$
\iota_\Lambda: 
X_\Lambda \longrightarrow
X_{{\frak L}^\Lambda},
$
which is the weakest topology for which $\iota_\Lambda$ is continuous.
Denote by 
$\widetilde{X}_\Lambda$
the topological space $X_\Lambda$ with the topology.
If $\Lambda$ is a topological Markov shift, the induced topology of 
$\widetilde{X}_\Lambda$
coincides with the original topology of 
$X_\Lambda$.
\proclaim{Lemma 7.3}
The topological space
$\widetilde{X}_\Lambda$ is generated by the clopen sets of the form
$U_\mu \cap \sigma_\Lambda^{-k}(\sigma_\Lambda^l(U_\nu))$
for
$\mu \in B_k(X_\Lambda), \nu \in B_l(X_\Lambda)$ with $k \le l$.
Hence the correspondence
$\chi_{U_\mu \cap \sigma_\Lambda^{-k}(\sigma_\Lambda^l(U_\nu))}
\longleftrightarrow S_\mu S_\nu^* S_\nu S_\mu^*$
yields an isomorphism
between
$C(\widetilde{X}_\Lambda)$ and ${\Cal D}_{{\frak L}^\Lambda}$.
\endproclaim
By the above lemma,
we know that
$C(\widetilde{X}_\Lambda)$ 
is isomorphic to
$C(X_{{\frak L}^\Lambda})$.

\smallskip

Let $\Lambda_1$ and 
$\Lambda_2$ be subshifts, and 
$X_{\Lambda_1}$ and 
$X_{\Lambda_2}$
their right one-sided subshifts.

\noindent
{\bf Definition.}
The subshifts
$(X_{\Lambda_1}, \sigma_{\Lambda_1})$ 
and 
$(X_{\Lambda_2},\sigma_{\Lambda_2})$
are said to be $\lambda$-{\it continuously orbit equivalent}\
if
there exists a homeomorphism
$h:X_{\Lambda_1}\longrightarrow X_{\Lambda_2}$,
that is also homeomorphic from
$\widetilde{X}_{\Lambda_1}\longrightarrow \widetilde{X}_{\Lambda_2}$
and
there exist continuous functions
$k_1,l_1:\widetilde{X}_{\Lambda_1}\longrightarrow \Zp$
and
$k_2,l_2:\widetilde{X}_{\Lambda_2}\longrightarrow \Zp$
such that
$$
\align
\sigma_{\Lambda_2}^{k_1(a)}(h\circ \sigma_{\Lambda_1}(a)) 
& =
\sigma_{\Lambda_2}^{l_1(a)}(h(a)) \quad\text{ for } a \in X_{\Lambda_1} \tag 7.1\\
\sigma_{\Lambda_1}^{k_2(b)}(h^{-1}\circ \sigma_{\Lambda_2}(b)) 
& =
\sigma_{\Lambda_1}^{l_2(b)}(h^{-1}(b)) \quad\text{ for } b \in X_{\Lambda_2}.\tag 7.2
\endalign
$$
We note that the conditions (7.1) and (7.2)
imply that 
$$
h(orb_{\sigma_{\Lambda_1}}(a)) = orb_{\sigma_{\Lambda_2}}(h(a)),
\quad
h^{-1}(orb_{\sigma_{\Lambda_2}}(b)) = orb_{\sigma_{\Lambda_1}}(h^{-1}(b))
$$
for
$a \in X_{\Lambda_1}, b \in X_{\Lambda_2}$.
\proclaim{Lemma 7.4}
Let
${\frak L}_1$ and
${\frak L}_2$
be the canonical $\lambda$-graph systems for 
$\Lambda_1$ and $\Lambda_2$ respectively. 
The following are equivalent:
\roster
\item
The subshifts 
$(X_{\Lambda_1},\sigma_{\Lambda_1})$
and
$(X_{\Lambda_2},\sigma_{\Lambda_2})$
are $\lambda$-continuously orbit equivalent.
\item
The factor maps
$\pi^{{\frak L}_1}_{\Lambda_1}$
and
$\pi^{{\frak L}_2}_{\Lambda_2}$
are continuously orbit equivalent.
\endroster
\endproclaim
\demo{Proof}
$(2)\Rightarrow (1)$ is clear.

$(1)\Rightarrow (2)$:
It suffices to show the equalities
$$
\align
\sigma_{{\frak L}_2}^{k_1(x)}(h(\sigma_{{\frak L}_1}(x))) 
& =
\sigma_{{\frak L}_2}^{l_1(x)}(h(x)), \quad\text{ for } x \in X_{{\frak L}_1} \\
\sigma_{{\frak L}_1}^{k_2(y)}(h^{-1}(\sigma_{{\frak L}_2}(y))) 
& =
\sigma_{{\frak L}_1}^{l_2(y)}(h^{-1}(y)), \quad\text{ for } y \in X_{{\frak L}_2}.
\endalign
$$
For $x \in X_{{\frak L}_1}$,
put
$k = k_1(x), l = l_1(x)$.
Since
$k_1,l_1:X_{{\frak L}_1} \longrightarrow \Zp$
are
continuous,
the set
$
U = \{z \in X_{{\frak L}_1} \mid k_1(z) = k, l_1(z) = l \}
$
is a clopen set in $X_{{\frak L}_1}$.
Since $X_{\Lambda_1}$ is dense in $X_{{\frak L}_1}$
through $\iota_{\Lambda_1}$,
one sees  $x \in U$ with $U \cap X_{\Lambda_1} \ne \emptyset$
and
the equality 
$$
\sigma_{{\frak L}_2}^{k_1(x)}(h\sigma_{{\frak L}_1}(x)) 
 =
\sigma_{{\frak L}_2}^{l_1(x)}(h(x)) \quad\text{ for } x \in X_{{\frak L}_1} 
$$
holds
because the equality holds for elements of $X_{\Lambda_1}$.
We similarly have the equality
$$
\sigma_{{\frak L}_1}^{k_2(y)}(h^{-1}\sigma_{{\frak L}_2}(y)) 
 =
\sigma_{{\frak L}_1}^{l_2(y)}(h^{-1}(y)) \quad\text{ for } y \in X_{{\frak L}_2}.
$$
Hence the factor maps
$\pi^{{\frak L}_1}_{\Lambda_1}$
and
$\pi^{{\frak L}_2}_{\Lambda_2}$
are continuously orbit equivalent.
\qed
\enddemo
Therefore we conclude:
\proclaim{Theorem 7.5}
Let $\Lambda_1$ and $\Lambda_2$ be subshifts satisfying condition (I).
The following are equivalent:
\roster
\item
There exists an isomorphism
$\Psi:\OLambdaA \longrightarrow \OLambdaB$
such that
$\Psi(\DLambdaA) = \DLambdaB$.
\item
The subshifts 
$(X_{\Lambda_1},\sigma_{\Lambda_1})$
and
$(X_{\Lambda_2},\sigma_{\Lambda_2})$
are $\lambda$-continuously orbit equivalent.
\endroster
\endproclaim
Let $A=[A(i,j)]_{i,j=1}^N$
be an $N \times N$
matrix with entries in $\{0,1\}$.
The Cuntz-Krieger algebra ${\Cal O}_A$ is generated by 
partial isometries $S_1,\dots,S_N$ 
satisfying 
$
\sum_{j=1}^N S_jS_j^* =1, S_i^*S_i = \sum_{j=1}^N A(i,j)S_jS_j^*, i=1,\dots,N.
$ 
The $C^*$-subalgebra generated by projections 
$
S_{\mu_n}^* \cdots S_{\mu_1}^* S_{\mu_1}\cdots S_{\mu_n},
\mu_1,\dots,\mu_n \in \{1,\dots,N\}
$
is canonically isomorphic to the commutative
$C^*$-algebra $C(X_A)$, that is denoted by  
$\DA$.
\proclaim{Corollary 7.6 ([Ma4])}
Let $A$ and $ B$ be square matrices with entries in $\{0,1\}$
 satisfying condition (I) in [CK]. 
Then the  following are equivalent:
 \roster
 \item There exists an isomorphism $\Psi: \OA \rightarrow  \OB$
 such that $\Psi(\DA) = \DB$.
\item
$(X_A, \sigma_A)$ and $(X_B, \sigma_B)$ are continuously orbit equivalent.
\endroster
\endproclaim
\demo{Proof}
For a topological Markov shift
$(X_A,\sigma_A)$,
the topology on
$\widetilde{X}_A$ coincides with the original topology on $X_A$.
Let $\Lambda_A$ be the two-sided topological Markov shift for the matrix $A$.
Then
$X_{\Lambda_A} = X_A$ 
and
${\Cal O}_{\Lambda_A} = \OA$
so that the assertion holds.
\qed
\enddemo 
Two one-sided subshifts 
$(X_{\Lambda_1},\sigma_{\Lambda_1})$
and
$(X_{\Lambda_2},\sigma_{\Lambda_2})$
are said to be topologically conjugate
if there exists a homeomorphism
$h:X_{\Lambda_1}\longrightarrow X_{\Lambda_2}$
such that
$\sigma_{\Lambda_2}\circ h = h \circ \sigma_{\Lambda_1}$,
and the homeomorphism $h$ is called a topological conjugacy.
One can prove that  
topological conjugacy gives rise to a $\lambda$-continuous orbit equivalence.
Hence we have.
\proclaim{Corollary 7.7([Ma3])}
Suppose that both subshifts 
 $\Lambda_1$ and $\Lambda_2$  satisfy condition (I).
Let 
$
h:(X_{\Lambda_1},\sigma_{\Lambda_1})
\rightarrow
(X_{\Lambda_2},\sigma_{\Lambda_2})
$
be a topological conjugacy of one-sided subshifts.
Then there exists an isomorphism
$\Psi:\OLambdaA\rightarrow \OLambdaB$
such that
$\Psi(\DLambdaA) = \DLambdaB$.
\endproclaim


\bigskip

{\it e-mail}: kengo{\@}yokohama-cu.ac.jp


\Refs
\refstyle{A}
\widestnumber\key{DGSW}


\ref
\no [Boy]
\by M. Boyle
\paper Topological orbit equivalence and factor maps in symbolic dynamics
\jour Ph. D. Thesis, University of Washington
\yr 1983
\endref

\ref
\no [BT]
\by M. Boyle and J. Tomiyama
\paper Bounded continuous orbit equivalence and $C^*$-algebras
\jour J. Math. Soc. Japan
\vol 50
\yr 1998
\pages 317--329
\endref


\ref
\no [CM]
\by T. M. Carlsen and K. Matsumoto
\paper Some remarks on the $C^*$-algebras associated with subshifts 
\jour Math. Scand.
\vol 95 
\yr 2004 
\pages 145--160
\endref


\ref
\no [CoKr]
\by A. Connes and W. Krieger
\paper Measure space automorphisms, the normalizers of their full groups,
and approximate finiteness 
\jour J. Funct. Anal. 
\vol 18
\yr 1975
\pages 318--327
\endref
 


\ref
\no [Cu]
\by J. Cuntz
\paper Simple $C^*$-algebras generated by isometries
\jour Comm. Math. Phys.
\vol 57
\yr 1977
\pages 173--185
\endref

\ref
\no [Cu2]
\by J. Cuntz
\paper Automorphisms of certain simple $C^*$-algebras
\jour in Quantum Fields-Algebras, Processes, Springer Verlag, Wien-New York 
\yr 1980
\pages 187--196
\endref 



\ref
\no [Cu3]
\by J. Cuntz 
\paper A class of $C^*$-algebras and topological Markov chains II: reducible chains and the Ext- functor for $C^*$-algebras
\jour Invent. Math.
\vol 63
\yr 1980
\pages 25--40
\endref



\ref
\no [CK]
\by J. Cuntz and W. Krieger
\paper A class of $C^*$-algebras and topological Markov chains
\jour Invent. Math.
\vol 56
\yr 1980
\pages 251--268
\endref


\ref 
\no [De]
\by V. Deaconu 
\paper Groupoids associated with endomorphisms
\jour Trans. Amer. Math. Soc. 
\vol 347
\yr 1995
\pages 1779--1786
\endref


\ref
\no [De2]
\by V. Deaconu 
\paper Generalized Cuntz-Krieger algebras
\jour Proc. Amer. Math. Soc. 
\vol 124
\yr 1996
\pages 3427--3435
\endref



\ref
\no [De3]
\by V. Deaconu 
\paper Generalized solenoids and $C^*$-algebras
\jour Pacific J. Math. 
\vol 190
\yr 1999
\pages 247--260
\endref



\ref
\no [De4]
\by V. Deaconu 
\paper Continuous graphs  and $C^*$-algebras
\jour Operator Theoretical Methods (Timi{\c c} soara, 1998) Theta Found.,
Bucharest\yr 2000
\pages 137--149
\endref




\ref
\no [D]
\by  H. Dye
\paper On groups of measure preserving transformations 
\jour American  J. Math.
\vol 81
\yr 1959
\pages 119--159
\endref

\ref
\no [D2]
\by  H. Dye
\paper On groups of measure preserving transformations II 
\jour American  J. Math.
\vol 85
\yr 1963
\pages 551--576
\endref






\ref
\no [GPS]
\by T. Giordano, I. F. Putnam and C. F. Skau
\paper Topological orbit equivalence and $C^*$-crossed products
\jour J. reine angew. Math.
\vol 469
\yr 1995
\pages 51--111
\endref


\ref
\no [GPS2]
\by T. Giordano, I. F. Putnam and C. F. Skau
\paper Full groups of Cantor minimal systems
\jour Isr. J. Math.
\vol 111
\yr 1999
\pages 285--320
\endref

\ref
\no [GMPS]
\by T. Giordano, H. Matui, I. F. Putnam and C. F. Skau
\paper Orbit equivalemce for Cantor minimal ${\Bbb Z}^2$-systems
\jour preprint
\yr 2006
\endref


\ref
\no [HO]
\by T. Hamachi and M. Oshikawa
\paper Fundamental homomorphisms of normalizer of ergodic transformation 
\jour Lecture Notes in Math. Springer
\vol 729
\yr 1978
\endref

\ref
\no [HPS]
\by R. H. Herman, I. F. Putnam and C. F. Skau
\paper Ordered Bratteli diagrams, dimension groups and topological dynamics 
\jour Internat. J. Math.
\vol 3
\yr 1992
\pages 827--864
\endref


\ref
\no [Ki]
\by  B. P. Kitchens
\book Symbolic dynamics
\publ Springer-Verlag
\publaddr Berlin, Heidelberg and New York
\yr 1998
\endref


\ref
\no [Kr]
\by W. Krieger
\paper On ergodic flows and isomorphisms of factors
\jour Math. Ann
\vol 223
\yr 1976
\pages 19--70
\endref





\ref
\no [Kr2]
\by W. Krieger
\paper On subshifts and topological Markov chains
\jour  Numbers, information and complexity (Bielefeld 1998),
Kluwer Acad. Publ. Boston MA (2000)
\pages 453--472
\endref


\ref
\no [KM]
\by W. Krieger and K. Matsumoto
\paper Shannon graphs, subshifts and lambda-graph systems
\jour  J. Math. Soc. Japan
\vol 54
\yr 2002
\pages 877--899
\endref


\ref
\no [LM]
\by  D. Lind and B. Marcus
\book An introduction to symbolic dynamics and coding
\publ Cambridge University Press
\publaddr Cambridge
\yr 1995
\endref



\ref
\no [Ma]
\by K. Matsumoto
\paper On $C^*$-algebras associated with subshifts
\jour  Internat. J. Math.
\vol 8
\yr 1997
\pages 357--374
\endref

\ref
\no [Ma2]
\by K. Matsumoto
\paper Presentations of subshifts and their topological conjugacy invariants
\jour Doc. Math.
\vol 4
\yr 1999
\pages 285-340
\endref


\ref
\no [Ma3]
\by K. Matsumoto
\paper On automorphisms of $C^*$-algebras associated with subshifts
\jour  J. Operator Theory
\vol 44
\yr 2000
\pages 91--112
\endref


\ref
\no [Ma4]
\by K. Matsumoto
\paper  $C^*$-algebras associated with presentations of subshifts
\jour  Doc. Math.
\vol 7
\yr 2002
\pages 1--30
\endref


\ref
\no [Ma5]
\by K. Matsumoto
\paper Orbit equivalence of topological Markov shifts and Cuntz-Krieger algebras\jour preprint, math  arXiv:0707.2114
\endref


\ref
\no [Ma6]
\by K. Matsumoto
\paper in preparation.
\endref

\ref
\no [Pat]
\by  A. L. T. Paterson
\book Groupoids, inverse semigroups, and their operator algebras
\publ Progress in Mathematics 170, Birkh{\"a}user
\publaddr Boston, Basel, Berlin
\yr 1998
\endref


\ref
\no [Put]
\by I. F. Putnam
\paper the $C^*$-algebras associated with minimal homeomorphisms of the Cantor set
\jour Pacific. J. Math. 
\vol 136
\yr 1989
\pages 329--353
\endref


\ref
\no [R{\o}1]
\by M. R{\o}rdam
\paper Classification of Cuntz-Krieger algebras 
\jour K-theory
\vol 9
\yr 1995
\pages 31--58
\endref


\ref
\no [To]
\by J. Tomiyama
\paper Topological full groups and structure of normalizers in transformation group $C^*$-algebras
\jour Pacific. J. Math. 
\vol 173
\yr 1996
\pages 571--583
\endref



\ref
\no [To2]
\by J. Tomiyama
\paper Representation of topological dynamical systems and $C^*$-algebras
\jour Contemporary  Math. 
\vol 228
\yr 1998
\pages 351--364
\endref
\bye